\documentclass[a4paper,fleqn]{cas-sc}

\usepackage{graphicx}
\usepackage{calc}
\usepackage{etoolbox}
\usepackage{xcolor}
\usepackage{calligra}

\makeatletter
\newlength{\hb@colwidth}
\def\makeheaderbox{
	\noindent\textcolor{black}{\rule{\textwidth}{1pt}}
	\par\vskip2pt\noindent
	\begin{minipage}{\textwidth}
		\centering
		\setlength{\tabcolsep}{0pt}
		\setlength{\extrarowheight}{0pt}
		\setlength{\hb@colwidth}{0.16666667\textwidth}
		\renewcommand{\arraystretch}{2}
		\begin{tabular}{
				p{\hb@colwidth}
				>{\columncolor[gray]{0.95}}p{4\hb@colwidth}
				p{\hb@colwidth}
			}
			\parbox[c][\hb@colwidth][c]{\hb@colwidth}{
				\centering
				\includegraphics[height=0.9\hb@colwidth, keepaspectratio]{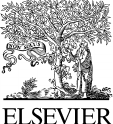}
			} & 
			\parbox[c][\hb@colwidth][c]{4\hb@colwidth}{
				\centering
				{\fontsize{8}{10}\selectfont Contents lists available at \href{https://www.elsevier.com/locate/fi}{ScienceDirect}}\par
				\vspace{15pt}
				{\fontsize{16}{18}\selectfont Journal of the Franklin Institute}\par
				\vspace{12pt}
				{\fontsize{8}{10}\selectfont journal homepage: \href{https://www.elsevier.com/locate/fi}{www.elsevier.com/locate/fi}}
			} & 
			\parbox[c][\hb@colwidth][c]{\hb@colwidth}{
				\centering
				\includegraphics[height=\hb@colwidth, keepaspectratio]{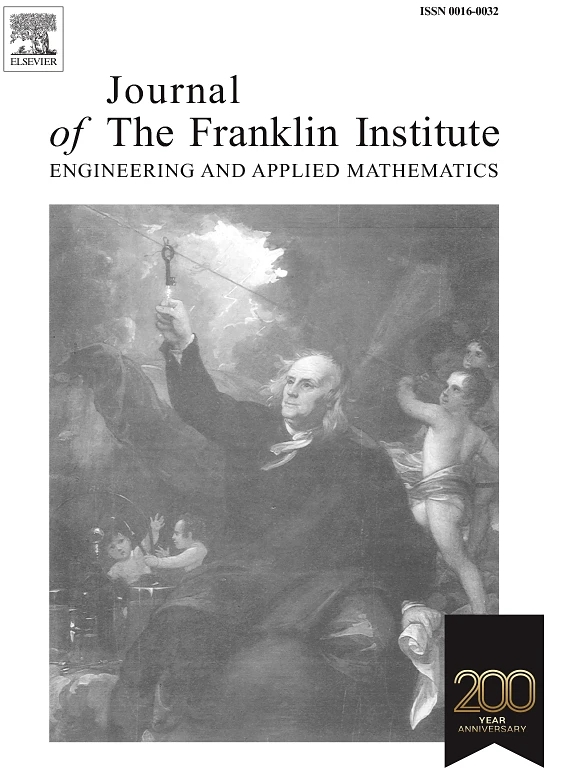}
			} \\
		\end{tabular}
	\end{minipage}
	\par\vskip2pt
	\noindent\textcolor{black}{\rule{\textwidth}{3pt}}
	\par\vskip18pt
}

\let\originalmaketitle\maketitle
\def\maketitle{
	\makeheaderbox
	\vspace{2pt}
	\begin{center}{\Large\bfseries \emph{Bicentennial Special Issue Invited Paper}}
    \end{center}\par
	\vspace{10pt}
	\originalmaketitle
}
\makeatother

\usepackage[numbers]{natbib}

\newtheorem{assum}{Assumption}

\def\tsc#1{\csdef{#1}{\textsc{\lowercase{#1}}\xspace}}
\tsc{WGM}
\tsc{QE}
\tsc{EP}
\tsc{PMS}
\tsc{BEC}
\tsc{DE}
\newcommand{\sgn}{\operatorname{sgn}}

\begin{document}
\let\WriteBookmarks\relax
\def\floatpagepagefraction{1}
\def\textpagefraction{.001}
\shorttitle{Journal of the Franklin Institute}
\shortauthors{T. R. Oliveira}





\title [mode = title]{Extremum Seeking Control:\\ Three Revolutions and the Road Ahead}

\tnotemark[1]

\tnotetext[1]{This document is the result of the research
   project funded by the Brazilian Funding Agencies: CNPq, CAPES, and FAPERJ.}


\author[1]{Tiago Roux Oliveira}[type=editor,
                        auid=000,bioid=1,
                        orcid=0000-0002-2232-8715
                       ]
\cormark[1]


\affiliation[1]{organization={Department of Electronics and Telecommunication Engineering, State University of Rio de Janeiro},
                addressline={RJ 20550-900}, 
                city={Rio de Janeiro},
                country={Brazil}}





%

\cortext[cor1]{Corresponding author}


\begin{abstract}
The history of extremum seeking is not merely the history of an algorithm; it is the history of an idea. Few ideas in control engineering have demonstrated the remarkable longevity of extremum seeking control (ESC). Invented more than one hundred years ago, ESC has continually reinvented itself while remaining faithful to its original objective: enabling systems to optimize their performance without relying on accurate mathematical models. This perspective article proposes that the evolution of ESC is best understood through three scientific revolutions. The first \textbf{Engineering Revolution (1922–1999)} established the engineering principles of model-free optimization; the second \textbf{Mathematical Revolution (2000–2010)} provided the rigorous mathematical foundations that transformed ESC into a mature discipline of nonlinear control; and the third ongoing \textbf{Infinite-Dimensional and Cyber-Physical Revolution (2010–present)} continues to expand its scope toward delays, partial differential equations, distributed optimization, event-triggered implementations, and increasingly complex cyber-physical systems. Beyond recounting this historical evolution, we offer a personal perspective on why ESC has remained relevant across successive technological eras. We argue that its enduring influence arises because the fundamental engineering challenge has never changed: ``how can a dynamical system learn to improve its own performance when the optimum is unknown?'' As optimization, learning, and feedback control become increasingly intertwined, ESC appears uniquely positioned to contribute to a new generation of intelligent autonomous systems, pointing toward what may become the field's \textbf{fourth scientific revolution}.
\end{abstract}




\maketitle

\section{Historical Perspectives: From Engineering Intuition to Learning Through Feedback}


\noindent \textit{It Really Is That Simple.}

Imagine you are cooking a pot of soup and trying to get the seasoning just right. You add a little salt, taste it, realize it still lacks flavor, add a little more, and taste it again—repeating this process until you reach the perfect balance.

This kind of trial and error comes naturally to us. But how can we enable a machine or an autonomous system to discover, entirely on its own, its \textbf{best operating point?}

This simple question has inspired more than a century of research in \textbf{extremum seeking control (ESC)}. ESC enables dynamical systems to adjust their behavior automatically so as to optimize their performance, even when the mathematical model describing the system is unknown or too complex to be used directly.

Rather than relying on an accurate model, the controller continuously experiments with small adjustments, observes their effect on performance, and gradually learns which direction leads to improvement. In this way, optimization becomes an integral part of the feedback loop itself.\\

\noindent \textit{The Pursuit of Optimal Performance.} 

Optimization is everywhere. Whether we seek to maximize energy efficiency, minimize fuel consumption, improve wireless communication, or increase industrial productivity, the underlying objective is always the same: to achieve the best possible performance under the prevailing conditions.

In practice, however, this goal is rarely straightforward. Many engineering systems operate under conditions that are only partially understood, continuously changing, or simply too complex to model with sufficient accuracy. Small variations in temperature, pressure, material properties, or environmental conditions may significantly alter system behavior. Developing an accurate mathematical model is often difficult—and sometimes practically impossible.

Extremum seeking control offers an elegant alternative. Rather than relying on a detailed mathematical description of the system, it learns directly from the system itself. By introducing small probing signals and observing the resulting changes in performance, the controller gradually identifies the direction in which performance improves, continuously steering the system toward its optimum.

This model-free philosophy is one of the defining characteristics of extremum seeking control and one of the main reasons for its enduring appeal across such a broad range of applications. \\

\noindent \textit{How Extremum Seeking Control Works: A Pinch of Intuition---and a Touch of Mathematics.} 

The basic idea behind extremum seeking is remarkably intuitive. Imagine a cook repeatedly tasting a soup while adjusting its seasoning. After each small modification, the cook evaluates whether the flavor has improved or deteriorated, using this information to decide the next adjustment. Through this simple process of experimentation and feedback, the recipe gradually converges to the desired taste.

Extremum seeking follows exactly the same principle. The controller deliberately introduces a small periodic perturbation into the control input and monitors the resulting variation in system performance. If the perturbation produces an improvement, the controller continues moving in that direction; if performance deteriorates, it adjusts in the opposite direction. Repeating this process continuously allows the system to approach the optimal operating point without ever knowing where that optimum lies beforehand.

From a mathematical viewpoint, the periodic perturbation plays a much deeper role than merely exploring the neighborhood of the current operating point. It acts as a real-time probe from which gradient information can be extracted through appropriate filtering and demodulation. In this way, the controller constructs an estimate of the local gradient of the unknown objective function and uses it to drive the optimization process.

Unlike many machine-learning algorithms, which often require extensive offline training and large datasets, extremum seeking learns while the system is operating. Optimization and control occur simultaneously, allowing the controller to continuously adapt to changes in the environment without interrupting normal operation.

This capability makes extremum seeking particularly attractive for systems operating under uncertain or time-varying conditions, including autonomous vehicles, industrial processes, communication networks, renewable energy systems, and many other modern engineering applications.

\begin{figure}
	\centering
	\includegraphics[width=.5\textwidth]{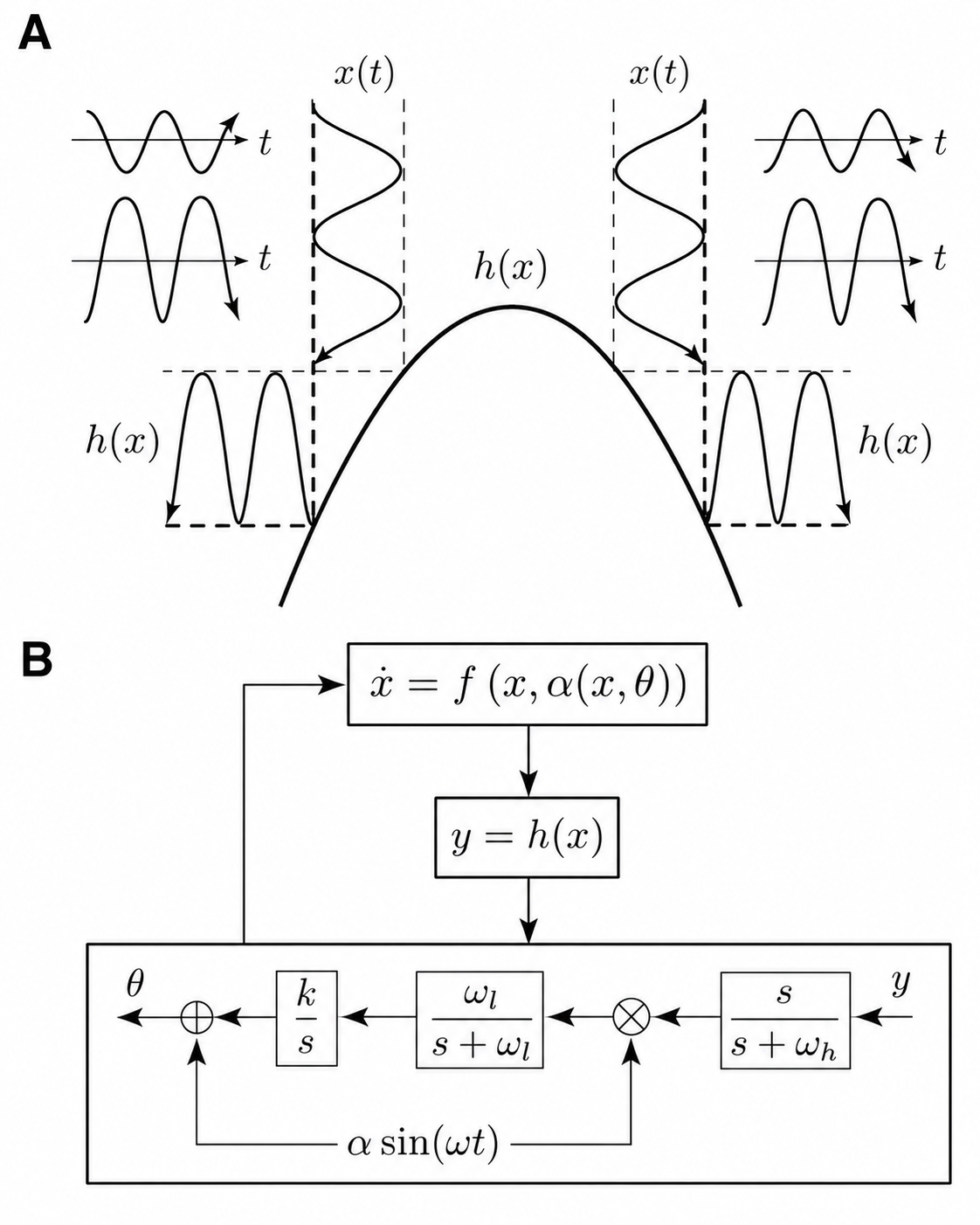}
	\caption{An intuitive sketch of the mechanism in searching for the maximum of
$h(x)$ near an extremum at which $h(x) \approx -x^2$ is also shown (A). Depending on whether $x(t) > x^*$ or $x(t) < x^*$ 
the value of the perturbed function $h(x(t))$ will
be in or out of phase with the perturbing signal. The ESC setup from \cite{KrsticWang2000} is shown for optimization of general nonlinear dynamic systems (B). Adapted from \cite{Scheinker2024}.}
	\label{figure1}
\end{figure}

%
Figure~\ref{figure1}A illustrates the underlying optimization mechanism and provides a simple intuitive explanation of why LeBlanc's method works. Starting with some initial condition $x(0)$, the system's input is perturbed in the form $x(t)=x(0)+\sin(\omega t)$, which results in a perturbation of some unknown output function, $h(x(t))=h(x(0)+\sin(\omega t))$. If $h(x)$ has an extremum point near $x(0)$ then $h(x(t))$ will be in or out of phase with the perturbing term $\sin(\omega t)$, as shown in Figure~\ref{figure1}A. The phase of the output function, relative to the perturbing term then gives the controller an idea of which way to move the parameter in order to approach the extremum point. 
Hence, the sinusoidal perturbation serves simultaneously as an exploration signal and as the key ingredient for estimating the local gradient of the objective function. Depending on whether the operating point lies to the left or to the right of the optimum, the resulting feedback drives the system toward the maximum. Figure~\ref{figure1}B shows the general extremum-seeking architecture for nonlinear dynamical systems---see also Appendix~\ref{fundamentals}. \\

\newpage
\noindent \textit{A Century of Extremum Seeking: From Engineering Intuition to Mathematical Theory.} 

The origins of extremum seeking control date back more than one hundred years. Surprisingly, one of the earliest known implementations was proposed as early as \textbf{1922}, when the French engineer \textbf{Maurice LeBlanc} introduced an automatic optimization scheme for railway electrification systems \cite{Leblanc1922}. Although remarkably simple, the algorithm already contained the essential ingredients that continue to define extremum seeking today: an integrator combined with a sinusoidal perturbation used simultaneously for exploration and gradient estimation.

At the time, the method was known as \textbf{peak seeking}, a name that explicitly reflected its objective—locating the operating point that maximized system performance. Despite its simplicity, the underlying idea proved remarkably powerful: by continuously perturbing the input and observing the resulting output variation, a controller could improve performance without requiring an explicit mathematical model of the process being optimized.

The concept attracted renewed attention in the early 1950s, when Draper and Li introduced peak-seeking methods to the American control community \cite{DraperLi1951}. Similar ideas had already appeared in the Soviet literature several years earlier, highlighting an interesting example of parallel scientific development during the early years of automatic control.

Throughout the 1950s and 1960s, numerous variants of extremum-seeking algorithms were proposed. Most were ingenious from an engineering standpoint and demonstrated impressive practical performance. Yet they largely remained heuristic. There was little theoretical understanding of why these algorithms worked, under what conditions they converged, or how their performance could be rigorously characterized.

For several decades, extremum seeking remained something of an engineering curiosity: a technique that practitioners knew often worked remarkably well, but whose mathematical foundations remained elusive.

Everything changed around the turn of the century. \\

\newpage
\noindent \textit{The Mathematical Revolution.} 

In 1995, {\AA}str{\"o}m and Wittenmark \cite{AstromWittenmark1995} identified extremum seeking as one of the most promising areas of adaptive control. A decisive breakthrough came with the pioneering work of Miroslav Krstic and Hsin-Hsiung Wang, whose landmark 2000 Automatica paper \cite{KrsticWang2000} established, for the first time, a rigorous nonlinear stability framework for extremum seeking. Their analysis demonstrated that the seemingly heuristic algorithm could, in fact, be understood within the well-established language of nonlinear systems theory.

The key insight was to recognize that extremum seeking naturally operates on two distinct time scales. The plant evolves according to its intrinsic dynamics, while the optimization process develops more gradually through the accumulation of information extracted from the periodic perturbation. By combining singular perturbation theory with averaging analysis, Krstic and Wang showed that the closed-loop system behaves, on average, like a gradient-based optimization algorithm.

This seemingly simple observation transformed the field. Instead of relying solely on empirical evidence, researchers could now establish rigorous guarantees of stability and convergence under precise mathematical assumptions. The analysis proved that the system converges exponentially to a neighborhood of the unknown optimum, whose size can be made arbitrarily small through appropriate choices of the perturbation parameters.

With this theoretical foundation in place, extremum seeking evolved from an ingenious engineering heuristic into a mature branch of nonlinear adaptive control. \\

\noindent \textit{From a Niche Technique to an International Research Field.} 

The impact of this theoretical breakthrough was immediate and far-reaching. What had once been a relatively specialized topic rapidly became one of the most active research areas at the intersection of adaptive control, optimization, and nonlinear systems.

During the early 2000s, the number of publications grew steadily. Within a decade, the community had expanded dramatically, with hundreds of new papers appearing every year. Today, the literature comprises \textbf{well over twenty thousand publications}, spanning virtually every branch of control engineering. This extraordinary evolution is comprehensively documented in Alex Scheinker's recent survey \cite{Scheinker2024}, \textit{100 Years of Extremum Seeking}, which provides an authoritative overview of the field's historical development, theoretical advances, and engineering applications.

Few areas of control have experienced such sustained growth over the past twenty-five years. The success of extremum seeking stems from a rare combination of characteristics: conceptual simplicity, rigorous mathematical foundations, broad applicability, and the ability to solve optimization problems online without requiring accurate system models \cite{editorial_CSm}. \\

\noindent \textit{Three Scientific Revolutions in Extremum Seeking.} 

Looking back over the past century, the evolution of extremum seeking can be understood as a sequence of three major scientific revolutions. The \textbf{first revolution (1922--1999)} was the birth of the engineering idea itself. Beginning with LeBlanc's pioneering work and continuing through several decades of practical developments, extremum seeking established itself as an ingenious model-free optimization strategy whose success was driven primarily by engineering intuition. The \textbf{second revolution (2000--2010)} occurred around the turn of the millennium, when nonlinear systems theory provided the rigorous mathematical foundations that transformed extremum seeking from an empirical technique into a mature discipline. Stability analysis, averaging theory, and singular perturbation methods made it possible to characterize the algorithm's behavior with mathematical precision. Today, we are witnessing\footnote{Throughout this Perspective, the plural ``we'' is used inclusively, inviting the reader to join the author in reflecting on the evolution and future of extremum seeking control; the singular ``I'' is reserved for personal experiences and reflections.} what may be regarded as the \textbf{third revolution (2010--present)}. The original framework, developed primarily for static maps and finite-dimensional dynamical systems, is being extended to increasingly complex environments involving delays, partial differential equations, distributed optimization, multi-agent systems, event-triggered implementations, stochastic dynamics, and learning-based architectures \cite{CSm_PDE,CSm_hybrid,CSm_industry,CSm_ML}. In many ways, the field is no longer concerned merely with finding optima—it is redefining what real-time optimization means for modern cyber-physical systems.  \\

\noindent \textit{The Third Revolution: Beyond Finite-Dimensional Systems.} 

Since the early 2000s, extremum seeking has expanded far beyond the framework for which the original theory was developed. Early algorithms focused primarily on static maps and finite-dimensional dynamical systems. As increasingly complex engineering problems emerged, however, the need arose to optimize systems whose dynamics are governed by delays, spatially distributed processes, communication constraints, and interactions among multiple decision makers.

These challenges required much more than incremental modifications of the classical algorithms. Delays fundamentally alter the information available to the optimizer, while partial differential equations introduce infinitely many dynamical modes. Distributed optimization brings together multiple agents with potentially competing objectives, each possessing only partial information about the global system. In all these settings, preserving the attractive model-free nature of extremum seeking while maintaining rigorous guarantees of stability and convergence becomes considerably more challenging.

Over the past two decades, significant progress has been achieved in addressing these problems. New theoretical tools have enabled extremum seeking to move beyond ordinary differential equations toward systems described by delay equations, transport phenomena, diffusion processes, wave propagation, and other infinite-dimensional dynamics. At the same time, the methodology has naturally evolved to encompass distributed optimization, non-cooperative games, event-triggered implementations, sampled-data architectures, stochastic optimization, and Newton-based acceleration techniques. \\

\noindent \textit{From Infinite-Dimensional Systems to Real Engineering Applications.} 

My own research has largely focused on this third stage in the evolution of extremum seeking. Rather than viewing delays and partial differential equations as obstacles to optimization, my work has sought to incorporate these dynamics directly into the controller design, thereby extending real-time optimization to systems previously considered beyond the reach of classical extremum-seeking methods.

This line of research culminated in the book \textit{Extremum Seeking through Delays and PDEs}, co-authored with Miroslav Krstic and published by SIAM in 2022 \cite{OliveiraKrstic2022}. The book presents a unified framework for extremum seeking in systems governed by delay differential equations and partial differential equations, including transport, heat, wave, and reaction-advection-diffusion models. It also extends the methodology from single-agent optimization to model-free Nash equilibrium seeking in heterogeneous multi-agent systems, demonstrating that optimization and game-theoretic learning can be performed even when the underlying dynamics are infinite-dimensional.

Perhaps more importantly, these theoretical developments have demonstrated that extremum seeking is not restricted to academic examples. The same mathematical principles can be applied across an extraordinarily diverse range of engineering problems.

My work has explored applications including optimal oil drilling, traffic-flow regulation, additive manufacturing, biological reactors, neuromuscular electrical stimulation, underwater source seeking using cable-actuated autonomous vehicles, light-source localization with flexible structures, and several other cyber-physical systems. Although these applications appear remarkably different, they all share a common objective: enabling complex dynamical systems to improve their performance autonomously while operating under uncertainty.

This universality is one of the most remarkable aspects of extremum seeking. Once viewed simply as an optimization algorithm, it has gradually evolved into a general framework for autonomous learning in dynamical systems.

\section{Pushing the Boundaries of Extremum Seeking}

Historical papers often reveal more than the state of knowledge of their time---they reveal the questions that researchers believed were worth asking. While preparing this article, one contribution particularly caught my attention: a paper published in the fifth issue of the \textit{IEEE Transactions on Automatic Control} in 1958 \cite{1104988}. Reading it today is a fascinating experience, not because of the algorithms it proposes, but because of the historical moment it captures.

At that time, adaptive control and nonlinear control were emerging aspirations rather than established disciplines. The now-familiar concepts of model-based adaptation, parameter estimation, and Lyapunov stability had not yet taken shape. Extremum seeking---then commonly referred to as the ``Draper–Li approach'' or simply as the ``hunt-and-seek-the-optimum'' method---represented one of the earliest visions of what adaptive control could become.

The author made an observation that, from today's perspective, is both insightful and surprisingly prophetic. He argued that model-free optimization through extremum seeking would become impractically slow for systems with many adjustable parameters and multiple time lags, suggesting that future adaptive controllers would inevitably require increasingly accurate mathematical models.

History unfolded differently.

\newpage
Model-based adaptive control indeed became one of the great successes of modern control theory, producing elegant mathematical frameworks and deep theoretical understanding. Yet, decades later, extremum seeking continues to find significant industrial application precisely because it requires so little prior knowledge of the system being optimized. If anything, the remarkable resurgence of extremum seeking reminds us that engineering history rarely follows a straight line. Ideas that appear limited in one era often return with renewed relevance when technology, computational power, and application demands evolve.\\

\noindent \textit{Yesterday's Limitations Often Become Tomorrow's Research Agenda.}

Perhaps the most intriguing aspect of this historical perspective is that some of the limitations identified in 1958 are no longer fundamental limitations today. Challenges associated with multiple inputs, communication delays, distributed decision-making, and complex dynamical systems have gradually become active research topics rather than barriers. In many respects, the recent evolution of extremum seeking can be viewed as a continuous effort to revisit those early questions using the mathematical and computational tools available to modern control theory.

The developments discussed in the following sections illustrate this evolution. They do not replace the original extremum-seeking paradigm; rather, they demonstrate how its underlying philosophy continues to adapt to increasingly challenging classes of dynamical systems.

The following examples illustrate how some of these long-standing challenges have gradually been addressed over the past decade. Rather than providing a comprehensive review, we focus on three representative research directions: \textbf{real-time optimization} itself, \textbf{Nash equilibrium seeking}, and \textbf{source seeking}. Further developments that, in our view, exemplify the ongoing expansion of the extremum-seeking paradigm beyond its classical formulation are also discussed later on.

\textcolor{black}{In order to develop ESC methods for systems with delays, the ESC
analysis was extended to systems governed by partial differential
equations (PDE), in which the delay is modeled as a \textit{transport hyperbolic} PDE \cite{Tiago}. Hence, the introduction to boundary control for ESC feedback is presented for hyperbolic-transport type of PDEs---see Appendix~\ref{SBA_Directory} (Section~\ref{backstepping_appendix}). On the other hand, there is no strong reason why the exposition provided could not have been conducted on some of the other classes of PDEs \cite{CSm_PDE,OLIVEIRA2023100908}. However, transport PDEs are particularly convenient because they are at the same time sufficiently simple and general to serve as a design template using which the reader can pursue extensions of the proposed designs and its integration to other classes of PDEs.} 

\subsection{Extremum Seeking through Delays}

In this section, we  derive the control algorithm for scalar gradient ESC in the presence of arbitrarily long input-output delays \cite{Tiago}. 
In our design we employ a predictor, to compensate for the delay. Since a predictor is always model-based, we are faced with the problem that we need to perform a prediction of a system in which either the input map or the output map contains a dependency on the Hessian, which is unknown. Hence, to perform a model-based prediction, \textit{i.e.}, delay compensation, we need an estimate of the Hessian. Our estimate of the Hessian is generated similarly to the estimation of the Hessian in Newton-based ESC \cite{GKN:2012}, namely, using a perturbation or demodulation signal. \\

\noindent \textit{Scalar Static Map with Delay.} 

As mentioned earlier, scalar ESC considers applications in which the goal is to
maximize (or minimize) the output $y\!\in\!\mathbb{R}$ of an unknown nonlinear static
map $Q(\theta)$ by varying the input $\theta\!\in\!\mathbb{R}$. Here, we additionally assume that there is a \textit{constant and known}
delay\index{constant delay} $D\geq0$ in the actuation path or measurement system such that the measured
output is given by
\begin{eqnarray} \label{delayed_output_old}
y(t)&=&Q(\theta(t-D))\,.
\end{eqnarray}
For notational clarity, we assume that our system is output-delayed in the following
presentation and block diagrams. However, the extension of the results in this article to the input-delay case is straightforward since, for a static map, any input delay can be moved to the output. The case when input delays $D_{\rm{in}}$ and output delays $D_{\rm{out}}$ occur
simultaneously could also be handled, by declaring that the total delay to be compensated is $D=D_{\rm{in}}+D_{\rm{out}}$, with $D_{\rm{in}}\,,D_{\rm{out}} \geq 0$.

%
Without loss of generality, let us consider the maximum seeking
problem such that the maximizing value of $\theta$ is denoted by
$\theta^*$. For the sake of simplicity, we also assume that the
nonlinear map is quadratic, \emph{i.e.},
\begin{eqnarray} \label{nonlinear_map}
Q(\theta)&=&y^*+\frac{H}{2}(\theta-\theta^*)^2\,,
\end{eqnarray}
where, besides the constants $\theta^*\!\in\!\mathbb{R}$ and $y^*\!\in\!\mathbb{R}$ being unknown,
the scalar $H\!<\!0$ is the unknown Hessian of the static map. 
By plugging (\ref{nonlinear_map}) into (\ref{delayed_output_old}), we obtain the \textit{quadratic static map with delay},
\begin{eqnarray} \label{delayed_output}
y(t)&=&y^*+\frac{H}{2}(\theta(t-D)-\theta^*)^2\,.
\end{eqnarray}

\begin{figure}
\begin{center}
\includegraphics[width=14cm]{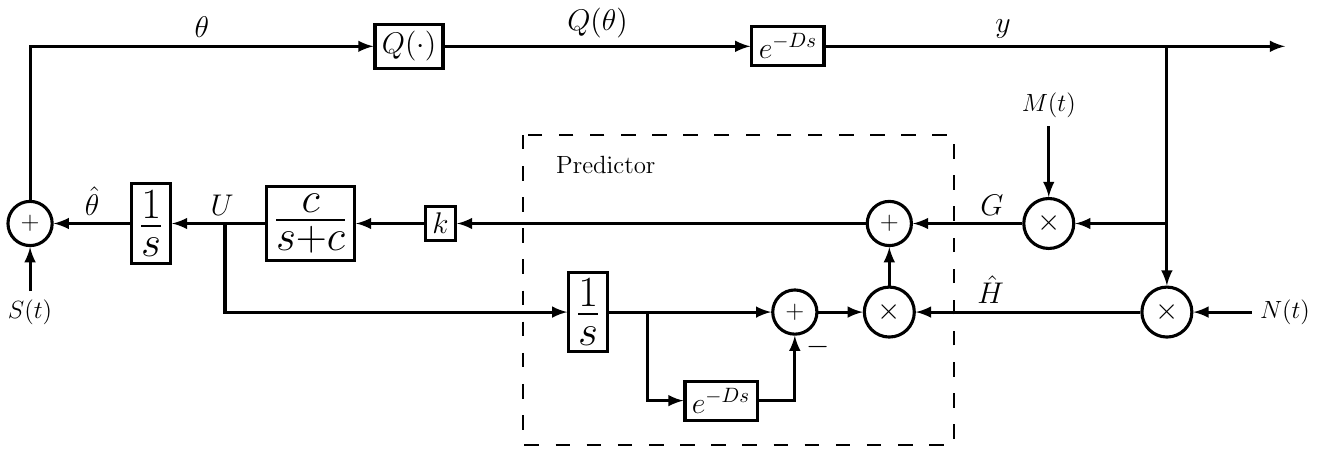}
\caption{Block diagram of the basic prediction scheme for
output-delay compensation in gradient ES. The predictor feedback with a perturbation-based estimate of the Hessian
obeys equation (\ref{predictor}), the probing and demodulation signals are given, respectively, by $S(t)=a \sin(\omega (t+D))$ and $M(t)=\frac{2}{a}
\sin(\omega t)$ and the demodulating signal $N(t)=-\frac{8}{a^2}\cos(2\omega t)$.}\label{Fig1}
\end{center}
\end{figure}

\noindent \textit{Probing and Demodulation Signals.} 

Let $\hat{\theta}$ be the estimate of $\theta^*$ and
\begin{eqnarray}\label{estimation_error}
\tilde{\theta}(t)=\hat{\theta}(t)-\theta^*
\end{eqnarray}
be the \textit{estimation error}. From Figure~\ref{Fig1} and Eq. (\ref{estimation_error}), the \textit{error dynamics} can be written as
\begin{eqnarray}
\dot{\tilde{\theta}}(t)=\dot{\hat{\theta}}(t)=U(t) \quad \mbox{and} \quad \dot{\tilde{\theta}}(t-D)=U(t-D)\,. \label{distinct_delayed_equation}
\end{eqnarray}
Moreover, one has
\begin{eqnarray}
G(t)&=& M(t)y(t)\,, \label{G} \\
\theta(t) &=& \hat{\theta}(t)+ S(t) \,, \label{theta}
\end{eqnarray}
where
the perturbation signal and demodulation signal, respectively, are given by
\begin{eqnarray}
S(t) &=& a \sin(\omega (t+D))\,, \label{dither_signalS}\\
M(t)&=& \frac{2}{a} \sin(\omega t) \label{dither_signalM}
\end{eqnarray}
with nonzero perturbation amplitude $a$ and frequency
$\omega$.

Notice that (\ref{dither_signalS}) is different from the
 perturbation signal used in the standard gradient extremum seeking free of delays.
Since the output delay can be transferred to the integrator
output for analysis purposes (or, equivalently, to its input), the phase shift
$+\omega D$ in (\ref{dither_signalS}) is applied to compensate the delay effect in the perturbation signal $S(t)$. The resulting additive dither signal $S(t)$ that anticipate propagation through delays are generally obtained by solving PDE motion-planning problems for more general PDE cases rather than transport equations, as discussed in \cite[Ch.~12]{krstic2008boundary}.

The signal
\begin{equation}\label{hessian_estimate_H}
\hat{H}(t)=N(t)y(t)\,,
\end{equation}
is used to obtain an estimate of the unknown Hessian $H$, where
the demodulating signal $N(t)$ is given by
\begin{equation}\label{hessian_estimate_N}
N(t)=-\frac{8}{a^2}\cos(2\omega t)\,.
\end{equation}
In \cite{GKN:2012}, it was proved that
\begin{equation}\label{averagingH}
\frac{1}{\Pi} \int_{0}^{\Pi}N(\sigma)y d\sigma = H\,, \quad \Pi=2\pi/\omega\,,
\end{equation}
if a quadratic map as in (\ref{nonlinear_map}) is considered. In
other words, the average version of (\ref{hessian_estimate_H}) is $\hat H_{\rm{av}}=(Ny)_{\rm{av}}=H$, \textit{i.e.}, the unknown Hessian is correctly estimated, on the average, using a product of available signals, the measured $y(t)$ and the given $N(t)$. \\

\noindent \textit{Predictor Feedback with a Perturbation-based Estimate of the Hessian}

By using the averaging analysis, we can verify that the average version of the signal (\ref{G}) is 
\begin{equation}\label{key_gradient}
G_{\rm{av}}(t)=H\tilde{\theta}_{\rm{av}}(t-D)\,.
\end{equation}
From (\ref{distinct_delayed_equation}) and (\ref{key_gradient}),
the following average models are also  obtained
\begin{eqnarray}\label{thetadot}
\dot{\tilde{\theta}}_{\rm{av}}(t-D)&=&U_{\rm{av}}(t-D)\,, \\
\dot G_{\rm{av}}(t)&=& H U_{\rm{av}}(t-D)\,, \label{Gdot}
\end{eqnarray}
where $U_{\rm{av}}\in \mathbb{R}$ is the resulting average control for $U\in \mathbb{R}$.

According to \cite{K:2009}, the delay in (\ref{thetadot}) can be represented using a transport PDE as (see Appendix~\ref{backstepping_appendix})
\begin{eqnarray}
\dot{\tilde{\theta}}(t-D) &=& u(0,t)\,, \label{PDE0_old}\\
u_t(x,t)&=&u_x(x,t)\,, \quad x\in[0,D]\,, \label{PDE1_old}\\
u(D,t)&=&U(t)\,, \label{PDE2_old}
\end{eqnarray}
where the solution of (\ref{PDE1_old})--(\ref{PDE2_old}) is
\begin{equation} \label{solution_nonaverage}
u(x,t) = U(t+x-D)\,.
\end{equation}
Here, we denote the partial derivatives of a function $u(x,t)$ as $\partial_x u(x,t) = \partial u(x,t)/\partial x,\ \partial_t u(x,t) = \partial u(x,t)/\partial t$. We conveniently use the compact form $u_x(x,t)$ and $u_t(x,t)$ for the former and the latter, respectively.

In order to motivate the predictor feedback design, the idea here is to compensate for the delay
by feeding back the future state $G(t+D)$, or $G_{\rm{av}}(t+D)$ in the equivalent average system.

Given any stabilizing gain $k>0$ for the undelayed system, our wish is to have a control that achieves
\begin{equation} \label{controlrewritten}
U_{\rm{av}}(t)=kG_{\rm{av}}(t+D)\,, \quad \forall t \geq 0\,,
\end{equation}
and this feedback law appears to be nonimplementable since it requires future values of the state. However, by applying the variation of constants
formula\index{variation of constants
formula} to (\ref{Gdot}), we can express the future state as 
\begin{equation}\label{prediction1}
G_{\rm{av}}(t+D)=G_{\rm{av}}(t)+H\int_{t}^{t+D}U_{\rm{av}}(\sigma-D)
d\sigma\,,
\end{equation}
where the current state $G_{\rm{av}}(t)$ is the initial condition.
Shifting the time variable under the integral in
(\ref{prediction1}), we obtain---for more details, see Appendix~\ref{predictor_appendix}:
\begin{equation}\label{prediction2}
G_{\rm{av}}(t+D)=G_{\rm{av}}(t)+H\int_{t-D}^{t}U_{\rm{av}}(\sigma) d\sigma\,,
\end{equation}
which gives the future state $G_{\rm{av}}(t+D)$ in terms of the
average control signal $U_{\rm{av}}(\sigma)$ from the past
window $[t-D,t]$. 
It yields the following feedback law\index{predictor}
\begin{equation}\label{predictoraverage}
U_{\rm{av}}(t)= k \left[G_{\rm{av}}(t)+H\int_{t-D}^{t}U_{\rm{av}}(\sigma) d\sigma
\right]\,.
\end{equation}

Hence, from (\ref{prediction2}) and (\ref{predictoraverage}),
the average feedback law (\ref{controlrewritten}) can be obtained indeed
as desired. Consequently,
\begin{equation} \label{hatthetarewritten}
\dot{\tilde{\theta}}_{\rm{av}}(t) = k G_{\rm{av}}(t+D)\,, \quad \forall t
\geq 0\,.
\end{equation}
%
Therefore, from (\ref{key_gradient}), one has
\begin{equation} \label{averaging_compensated}
\frac{d \tilde{\theta}_{\rm{av}}(t)}{d t} = k H
\tilde{\theta}_{\rm{av}}(t)\,, \quad \forall t \geq D\,,
\end{equation}
with an exponentially attractive equilibrium
$\tilde{\theta}_{\rm{av}}^{e}\!=\!0$, since $k\!>\!0$ in the
control design and $H<0$ by assumption. It means that the delay is perfectly compensated
in $D$ units of time, namely, the system evolves as if the delay were
absent after $D \geq 0$.

In the next section, we show that the control objectives can still be achieved if a simple modification
of the above basic predictor-based controller, which employs a low-pass filter, is applied. In this case,
we propose the following infinite-dimensional and averaging-based predictor feedback\index{filtered predictor} in order to compensate the delay \cite{K:2008}
\begin{eqnarray}\label{predictor}
U(t)&=& \frac{c}{s+c}\left\{k \left[G(t)+\hat{H}(t)\int_{t-D}^{t} U(\tau) d\tau \right] \right\}\,,
\end{eqnarray}
where $c>0$ is sufficiently large, \textit{i.e.},
the predictor feedback is a low-pass-filtered version of the non-averaged controller in 
%
(\ref{predictoraverage}). This low-pass filtering is particularly required in the stability analysis when the averaging theorem in infinite dimensions \cite{HL:1990,L:2002} is invoked. Note that we mix the time- and frequency-domain notation in (\ref{predictor}) by using the braces $\{\cdot\}$ to denote that the transfer function acts as an operator on a time-domain function.

The predictor feedback (\ref{predictor}) is infinite-dimensional because the integral involves
the control history over the interval $[t-D,t]$. This feedback is also averaging-based (perturbation-based) because $\hat{H}$ is
updated according to the estimate (\ref{hessian_estimate_H}) of the unknown Hessian $H$, which satisfies
the averaging property (\ref{averagingH}).

\subsection{Games with Delays}

\noindent \textit{$N$-Player Game with Quadratic Payoffs for Nash Equilibrium Seeking.} 

Game theory provides an important framework for mathematical modeling and analysis of scenarios involving different players where there is coupling in their actions, in the sense that their respective outcomes (outputs) $y_{i}(t) \in \mathbb{R}$ do not depend exclusively on their own actions/strategies (input signals) $\theta_{i}(t)\in \mathbb{R}$, with $i=1 ,\ldots, N$, but at least on a subset of others'. Moreover, defining $\theta := [\theta_1, \ldots, \theta_N]^T$, each player's payoff function $J_{i}(\theta) : \mathbb{R}^{N} \to \mathbb{R}$ depends on the action $\theta_j$ of at least one other Player $j$, $j\not= i$. An $N$-tuple of actions, $\theta^*$, is said to be in Nash equilibrium, if no Player $i$ can improve its payoff  by unilaterally deviating from $\theta_i^*$, this being so for all $i$ \cite{Basar:1999}. 

Consider games where the payoff function of each player is quadratic, 
expressed as a strictly concave\footnote{By strict concavity, we mean $J_i(\theta)$ is strictly concave in $\theta_i$ for all $\theta_{-i}$, this being so for each $i=1,\ldots, N$.} combination of their actions 
\begin{align}
J_{i}(\theta(t))=&\frac{1}{2}\sum_{j=1}^{N}\sum_{k=1}^{N}\epsilon_{jk}^{i}H_{jk}^{i}\theta_{j}(t)\theta_{k}(t) 
+\sum_{j=1}^{N}h_{j}^{i}\theta_{j}(t)+c_{i}\,, \label{eq:Ji}
\end{align}  
where $\theta_{j}(t) \!\in\! \mathbb{R}$ is the decision variable (action) of Player $j$, 
$H_{jk}^{i}$, $h_{j}^{i}$, $c_{i} \!\in\! \mathbb{R}$ are constants, $H_{ii}^{i}\!<\!0$, $H_{jk}^{i}\!=\!H_{kj}^{i}$ and $\epsilon_{jk}^{i}\!=\!\epsilon_{kj}^{i}\!>\!0$, $\forall i,j,k$. 

Quadratic payoff functions, of the type above, are of particular interest in game theory, first because they constitute second-order approximations to other types of non-quadratic payoff functions, and second because they are analytically tractable, leading in general to closed-form equilibrium solutions which provide insight into the properties and features of the equilibrium solution concept under consideration \cite{Basar:1999}. 

For the sake of completeness, we provide here in mathematical terms, the definition of a Nash equilibrium $\theta^*=[\theta^*_1\,, \ldots\,,\theta_N^*]^T$ in an $N$-player game:
\begin{eqnarray} \label{Nashcu}
J_i(\theta_i^*\,,\theta_{-i}^*) \!\geq\! J_i(\theta_i\,,\theta_{-i}^*)\,, \quad 
\forall \theta_i \in \mathcal{U}_i\,, \quad i \!\in\! \{1\,, \ldots\,,N\}\,,
\end{eqnarray}
where $J_i$ is the payoff function of the player $i$, the variable $\theta_i$ corresponds to its action, while $\mathcal{U}_i$ is its action set and $\theta_{-i}$ denotes the actions of the other players. Hence, no player has an incentive to unilaterally deviate its action from $\theta^*$. 

In a duopoly example ($N=2$), as considered in Appendix~\ref{games_appendix}, $\mathcal{U}_1=\mathcal{U}_2=\mathbb{R}$, where $\mathbb{R}$ denotes the set of real numbers. 
In order to determine the Nash equilibrium solution in strictly concave quadratic games with $N$ players, where each action set is the entire real line, one should differentiate $J_{i}$ with respect to $\theta_{i}(t) \,, \forall i=1 ,\ldots, N$, setting the resulting expressions equal to zero and solving the set of equations thus obtained.

This set of equations, which also provides a sufficient condition due to the strict concavity, is
\begin{align}
\sum_{j =1}^{N}\epsilon_{ij}^{i}H_{ij}^{i}\theta_{j}^{*}+h_{i}^{i}=0\,,\quad i=1 ,\ldots, N\,, \label{eq:NE_v0}
\end{align}
which can be written in compact form as 
\begin{align}
\begin{bmatrix}
\epsilon_{11}^{1} H_{11}^{1} & \epsilon_{12}^{1} H_{12}^{1} & \hdots & \epsilon_{1N}^{1} H_{1N}^{1} \\
\epsilon_{21}^{2} H_{21}^{2} & \epsilon_{22}^{2} H_{22}^{2} & \hdots & \epsilon_{2N}^{2} H_{2N}^{2} \\
\vdots                       & \vdots                       &        & \vdots     \\
\epsilon_{N1}^{N} H_{N1}^{N} & \epsilon_{N2}^{N} H_{N2}^{N} & \hdots & \epsilon_{NN}^{N} H_{NN}^{N}   
\end{bmatrix}
\begin{bmatrix}
\theta_{1}^{*} \\
\theta_{2}^{*} \\
\vdots \\
\theta_{N}^{*} 
\end{bmatrix}
=-
\begin{bmatrix}
h_{1}^{1} \\
h_{2}^{2} \\
\vdots    \\
h_{N}^{N}   
\end{bmatrix}
\,.
\end{align}
Defining the Hessian matrix $H$ and vectors $\theta^*$ and $h$ by 
%
\begin{align}
H&:=
\begin{bmatrix}
\epsilon_{11}^{1} H_{11}^{1} & \epsilon_{12}^{1} H_{12}^{1} & \hdots & \epsilon_{1N}^{1} H_{1N}^{1} \\
\epsilon_{21}^{2} H_{21}^{2} & \epsilon_{22}^{2} H_{22}^{2} & \hdots & \epsilon_{2N}^{2} H_{2N}^{2} \\
\vdots                       & \vdots                       &        & \vdots     \\
\epsilon_{N1}^{N} H_{N1}^{N} & \epsilon_{N2}^{N} H_{N2}^{N} & \hdots & \epsilon_{NN}^{N} H_{NN}^{N}   
\end{bmatrix}
\,, \nonumber \\
\theta^{*}&:=
\begin{bmatrix}
\theta_{1}^{*} \\
\theta_{2}^{*} \\
\vdots \\
\theta_{N}^{*} 
\end{bmatrix}
\,, \quad
h:=
\begin{bmatrix}
h_{1}^{1} \\
h_{2}^{2} \\
\vdots    \\
h_{N}^{N}   
\end{bmatrix}
\,, \label{eq:Htheta*h}
\end{align}
there exists a unique Nash Equilibrium at  $\theta^{*}=-H^{-1}h$, if $H$ is invertible:
\begin{align}
\begin{bmatrix}
\theta_{1}^{*} \\
\theta_{2}^{*} \\
\vdots \\
\theta_{N}^{*} 
\end{bmatrix}
\!=\!-
\begin{bmatrix}
\epsilon_{11}^{1} H_{11}^{1} & \epsilon_{12}^{1} H_{12}^{1} & \hdots & \epsilon_{1N}^{1} H_{1N}^{1} \\
\epsilon_{21}^{2} H_{21}^{2} & \epsilon_{22}^{2} H_{22}^{2} & \hdots & \epsilon_{2N}^{2} H_{2N}^{2} \\
\vdots                       & \vdots                       &        & \vdots     \\
\epsilon_{N1}^{N} H_{N1}^{N} & \epsilon_{N2}^{N} H_{N2}^{N} & \hdots & \epsilon_{NN}^{N} H_{NN}^{N}   
\end{bmatrix}^{-1}\!\!
\begin{bmatrix}
h_{1}^{1} \\
h_{2}^{2} \\
\vdots    \\
h_{N}^{N}   
\end{bmatrix}. \label{eq:NE_v2}
\end{align}
%
For more details, see \cite[Ch. 4]{Basar:1999}.

In addition to Assumption~\ref{ch.13.Assumption 1.} formulated in \cite{FKB:2012}, we further assume/formalize the following Assumption~\ref{ch13.Assumption 2.} for noncooperative games.

\begin{assum} \label{ch.13.Assumption 1.}
The Hessian matrix $H$ given by \textcolor{black}{\textnormal{(\ref{eq:Htheta*h})}} is strictly diagonal dominant, \textit{i.e.},
\begin{equation}
\sum_{j\neq i}^{N}|\epsilon_{ij}^{i}H_{ij}^{i}| < |\epsilon_{ii}^{i}H_{ii}^{i}|\,, \quad i \in \{1\,, \ldots N\}\,.    
\end{equation}
\end{assum}

\begin{assum} \label{ch13.Assumption 2.}
The parameters $\epsilon_{jk}^{i}$ and $\epsilon_{kj}^{i}$ which appear in the Hessian matrix $H$ given by (\ref{eq:Htheta*h}) satisfy the conditions below: 
\begin{eqnarray}
\epsilon_{ii}^{i}=1\,, \quad \epsilon_{jk}^{i}&=&\epsilon_{kj}^{i}=\epsilon\,, \quad \forall j\neq k\,,
\end{eqnarray}
with $0<\epsilon<1$ in the payoff functions  \textcolor{black}{\textnormal{(\ref{eq:Ji})}}.
\end{assum}

By Assumption~\ref{ch.13.Assumption 1.}, the Nash Equilibrium $\theta^*$ exists and is unique, since strictly diagonally dominant matrices are nonsingular by the Levy-Desplanques Theorem \cite{HJ:1985}. 
Assumption~\ref{ch13.Assumption 2.} could be relaxed to consider different values of the coupling parameters $\epsilon_i$ for each Player $i$. However, without loss of generality, we have assumed the same weights for the interconnection channels among the players in order to facilitate the proofs of our theorems, but also to guarantee that the considered noncooperative game is not favoring any specific player.\\

\noindent \textit{Noncooperative Scenario with Delays.}

\begin{figure}
\begin{center}
\includegraphics[width=3.0 in]{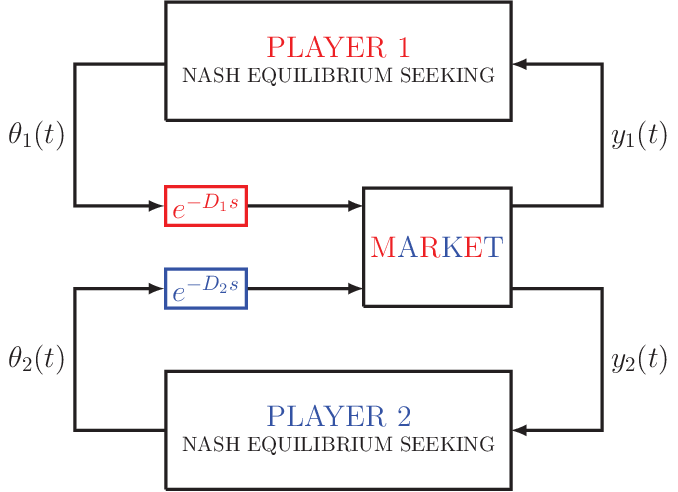}
\end{center}
\caption{Nash equilibrium-seeking schemes applied by two players ($N=2$) in a duopoly market structure with delayed players' actions.}
\label{fig:blockDiagram_v2}
\end{figure}

For the sake of completeness and to keep the material \textcolor{black}{self-contained}, we briefly review the case of noncooperative games subject to multiple and distinct delays originally addressed in \cite{JOTA:2021}. In this scenario, the purpose of the extremum seeking is still to estimate the Nash equilibrium vector $\theta^*$, but without allowing any sharing of information among the players. As mentioned earlier, each player only needs to measure the value of its own payoff function described by  
\begin{align} \label{ajagambiarra_noncooperative}
y_i(t)\!=\!J_i(\theta(t\!-\!D))=&\frac{1}{2}\sum_{j=1}^{N}\sum_{k=1}^{N}\epsilon_{jk}^{i}H_{jk}^{i}\theta_{j}(t-D_{j})\theta_{k}(t-D_{k}) +\sum_{j=1}^{N}h_{j}^{i}\theta_{j}(t-D_{j})+c_{i}\,,   
\end{align}
with $J_i$ given by (\ref{eq:Ji}). In this sense, we are able to formulate the closed-loop system in a \textit{decentralized} fashion, where no knowledge about the payoffs $y_{-i}$ or actions $\theta_{-i}$ of the other players is required, as illustrated in Figure~\ref{fig:blockDiagram_v2}. 
%
%

Without loss of generality, we assume that the inputs have distinct known (constant) delays which are ordered so that
\begin{equation}\label{equ:4_delay_order}
	D = \text{diag}\{D_1,D_2,\cdots, D_N\}, \quad 0\leq D_1 \leq \cdots \leq D_N\,.
\end{equation}
Further, given any $\mathbb{R}^N$-valued signal $f$, we introduce 
\begin{small}
\begin{equation} \label{4_fD}
	f^D(t) \!:=\! f(t-D)\!=\!  
	\begin{bmatrix}
		f_1(t\!-\!D_1) \!&\! f_2(t\!-\!D_2) \!&\! ... \!&\! f_N(t\!-\!D_N)
	\end{bmatrix}^T.
\end{equation} 
\end{small}

\begin{figure}
\begin{center}
\includegraphics[width=3.25 in]{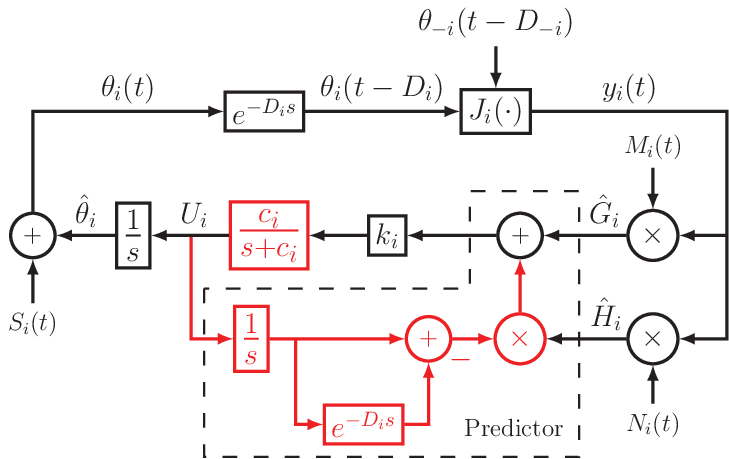}
\end{center}
\caption{Block diagram illustrating the Nash Equilibrium-seeking strategy performed for each player. \textcolor{black}{In red color, the predictor feedback used to compensate the individual delay $D_i$ for the noncooperative case. With some abuse of notation, the constants $c_i$ were chosen to denote the parameters of the filters $c_i/(s+c_i)$, but they are not necessarily the same constants which appear in the payoff functions given in (\ref{eq:Ji}).}}
\label{fig:blockDiagram_v3}
\end{figure}

Figure~\ref{fig:blockDiagram_v3} contains a schematic diagram that summarizes the proposed Nash Equilibrium policy for each $i$-th player where its output is given by (\ref{ajagambiarra_noncooperative}), 
where the vector $\theta_{-i}(t-D_{-i})$ in Figure~\ref{fig:blockDiagram_v3} represents the delayed actions of all other players.

The additive-multiplicative dithers $S_i(t)$, $M_i(t)$ are
\begin{align}
S_{i}(t)&=a_{i}\sin(\omega_{i}(t+D_{i})) \,, \label{eq:Si} \\
M_{i}(t)&=\frac{2}{a_{i}}\sin(\omega_{i}t)\,, \label{eq:Mi}
\end{align}
with nonzero constant amplitudes $a_i>0$ at frequencies $\omega_i \neq \omega_j$. 
%
Such probing frequencies $\omega_i$ can be selected as
\begin{equation}\label{omegadefinition}
\omega_i=\omega_{i}'\omega=\mathcal{O}(\omega)\,, \quad
i\in{1,2,\ldots, N}\,,
\end{equation}
where $\omega$ is a positive constant and $\omega_{i}'$ is a
rational number. One possible choice is given in \cite{GKN:2012} as
\begin{equation}\label{possible_implementation}
\omega_{i}'\not\in \left\{\omega_{j}'\,, \ \frac{1}{2}(\omega_{j}'+\omega_{k}')\,, \ \omega_{j}'+2\omega_{k}'\,, \ \omega_{j}'+\omega_{k}'\pm \omega_{l}'\right\}\,,
\end{equation}
for all distinct $i,j,k$ and $l$. 

Considering $\hat{\theta}_{i}(t)$ as an estimate of $\theta^{*}_{i}$, one can define the \textit{estimation error}:
\begin{align}
\tilde{\theta}_{i}(t)&=\hat{\theta}_{i}(t)-\theta_{i}^{*}\,.\label{eq:tildeThetai}
\end{align}
The estimate $\hat{G}_{i}$ of the unknown gradient of each payoff $J_i$ is given by
\begin{align}
\hat{G}_{i}(t)&=M_{i}(t)y_{i}(t)\,. \label{eq:hatGi_nonaverage}
\end{align}
%
Computing the average of the resulting signal leads us to
\begin{align}
\hat{G}_{i}^{\rm{av}}(t)
&=\sum_{j=1}^{N}\epsilon_{ij}^{i}H_{ij}^{i}\tilde{\theta}_{j}^{\rm{av}}(t-D_{j})\,,\label{eq:hatGi}
\end{align}
with $\Pi$ defined as
\begin{equation} \label{period}
\Pi:= 2 \pi \times \text{LCM}\left\{\frac{1}{\omega_i} \right\}\,, 
\end{equation}
and LCM standing for the least common multiple.

At this point, if we ignore the prediction loop and the low-pass filter (both indicated in red) in Figure~\ref{fig:blockDiagram_v3}, the control law $U_i(t)=k_i\hat{G}_i(t)$ could be obtained as in the classical ESC approach. In this case, from equations (\ref{eq:tildeThetai}) and (\ref{eq:hatGi}), we could write the average version of
\begin{equation} \label{cudoerro}
\dot{\tilde{\theta}}_i(t)=U_i(t)
\end{equation}
as 
\begin{align}
\dot{\tilde{\theta}}_{i}^{\rm{av}}(t)&=k_{i}\hat{G}_{i}^{\rm{av}}(t) \nonumber \\
&=k_{i}\sum_{j=1}^{N}\epsilon_{ij}^{i}H_{ij}^{i}\tilde{\theta}_{j}^{\rm{av}}(t-D_{j})\,.\label{dotTildeAvi_v0}
\end{align}
Therefore, by defining $$\tilde{\theta}^{\rm{av}}(t): = [\tilde{\theta}_{1}^{\rm{av}}(t)\,,\tilde{\theta}_{2}^{\rm{av}}(t)\,,\ldots\,,\tilde{\theta}_{N}^{\rm{av}}(t)]^T \in \mathbb{R}^{N}$$ and taking into account all players, one has
\begin{align}
\dot{\tilde{\theta}}^{\rm{av}}(t)&=KH\tilde{\theta}^{\rm{av}}(t-D)\,,\label{dotTildeAv_v0}
\end{align}
with $K:=\text{diag}\{k_{1}\,,\ldots\,,k_{N}\}$ and $H$ given by (\ref{eq:Htheta*h}). Equation (\ref{dotTildeAv_v0}) means that, even if $KH$ was a Hurwitz matrix, the equilibrium $\tilde{\theta}_{\rm{e}}^{\rm{av}}=0$ of the closed-loop average system would not necessarily be stable for arbitrary values of the time-delays $D_i$. This reinforces the demand to employ the prediction feedback $U_i(t)=k_i\hat{G}_i(t+D_i)$---or even its filtered version---for each player to collectively stabilize the closed-loop system, as illustrated with red color in Figure~\ref{fig:blockDiagram_v3}.

In such a decentralized scenario, the dither frequencies $\omega_{-i}$, the excitation amplitudes $a_{-i}$, and consequently the individual control laws $U_{-i}(t)$ are not available to Player $i$. Recalling that the model of the payoffs (\ref{eq:Ji}) and (\ref{ajagambiarra_noncooperative}) are also assumed to be unknown, it becomes impossible to reconstruct individually or estimate the Hessian matrix $H$ given in (\ref{eq:Htheta*h}) by using demodulating signals such as in \cite{Tiago}. 

%
%

Following the non-sharing information paradigm, the $i$th-player is only able to estimate the element $H_{ii}^{i}$ of the $H$ matrix (\ref{eq:Htheta*h}) by itself, and this being so for all players. Therefore, only the diagonal of $H$ can be properly recovered in the average sense. In this way, the signal $N_{i}(t)$ is now simply defined as:
\begin{align}
N_{i}(t):= N_{ii}(t)=	\dfrac{16}{a^2_i} \bigg( \sin^2 (\omega_i t) - \dfrac{1}{2} \bigg)\,, \label{eq:Ni}
\end{align} 
according to \cite{Tiago}. Then, the average version of 
\begin{align}
\hat{H}_{i}(t)= N_i(t) y_i(t) \label{eq:hatHi_nonaverage}
\end{align}
is given by
\begin{align}
\hat{H}_{i}^{\rm{av}}(t)=\left[N_{i}(t)y_{i}(t)\right]_{\rm{av}}=H_{ii}^{i}\,. \label{eq:hatHi}
\end{align}

In order to compensate for the time delays, the following predictor-based update law was proposed in \cite{JOTA:2021}:
\begin{align} \label{adaptation_NC}
	\dot{\hat{\theta}}_{i}(t) &= U_{i}(t)\,, \\
	\dot{U}_{i}(t) &= -c_{i}U_{i}(t) + c_{i}k_{i} \bigg(\hat{G}_{i}(t) +\hat{H}_{i}(t) \int_{t-D_i}^{t} U_i(\tau)d\tau \bigg)\,, \label{4_dU_NC}
\end{align}
\noindent for positive constants $k_{i}$ and $c_{i}$.

\subsection{Source Seeking with Delays}


\noindent \textit{Nonholonomic Mobile Robots}

Consider a unicycle-type mobile robot equipped with a sensor mounted on its structure. Unlike the classical ESC problem, the measurements provided by the sensor are assumed to be subject to a constant and known measurement delay $D > 0$. This delay represents the sensor response time, sampling latency, or processing time involved in signal acquisition.

The kinematic dynamics of the unicycle robot are described by
\begin{align}
    \dot{x}(t) &= v(t)\cos\theta(t) \label{eq:din_x}, \\
    \dot{y}(t) &= v(t)\sin\theta(t) \label{eq:din_y}, \\
    \dot{\theta}(t) &= \omega(t) \label{eq:din_theta},
\end{align}
where $x(t), y(t) \in \mathbb{R}$ represent the Cartesian coordinates of the center of the robot axle and $\theta(t) \in \mathbb{R}$ is its orientation angle, while $v(t)$ and $\omega(t)$ represent the linear and angular velocity control inputs, respectively.

\newpage
To overcome the nonholonomic constraint of the unicycle and enable control of the vehicle in the Cartesian plane, consider the position of a point of interest $p(t) = [x_p(t), y_p(t)]^\top \in \mathbb{R}^2$ located at a distance $l > 0$ ahead of the wheel center along the robot symmetry axis:
\begin{align}
    x_p(t) &= x(t) + l \cos\theta(t) \label{eq:ponto_xp}, \\
    y_p(t) &= y(t) + l \sin\theta(t) \label{eq:ponto_yp}.
\end{align}

Differentiating \eqref{eq:ponto_xp}--\eqref{eq:ponto_yp} with respect to time yields the kinematic relationship with the actual inputs $[v(t), \omega(t)]^\top$:
\begin{align}
    \begin{bmatrix} \dot{x}_p(t) \\[2pt] \dot{y}_p(t) \end{bmatrix} &= \begin{bmatrix} \cos\theta(t) & -l\sin\theta(t) \\ \sin\theta(t) & l\cos\theta(t) \end{bmatrix} \begin{bmatrix} v(t) \\[2pt] \omega(t) \end{bmatrix} = T(\theta(t)) \begin{bmatrix} v(t) \\[2pt] \omega(t) \end{bmatrix} \label{eq:matriz_decouplage},
\end{align}
where the decoupling matrix $T(\theta(t))$ is nonsingular for any $l > 0$.\\

\noindent \textit{Source Signal Modeling}

The objective of the robot is to locate a source emitting a signal whose intensity decreases with the distance from the optimal point in the plane $p^* = [x^*, y^*]^\top$. A local quadratic approximation of the scalar field in a neighborhood of the maximum is adopted:
\begin{align}
    Q(p(t)) = Q^{*} - \frac{q_x}{2}(x_p(t) - x^{*})^{2} - \frac{q_y}{2}(y_p(t) - y^{*})^{2} \label{eq:mapa_Q}, \qquad q_x>0, \quad q_y>0,
\end{align}
which has a unique local maximum $Q^* = Q(p^*)$.

Due to the intrinsic sensor measurement delay, the scalar reading available at time $t$ is given by $y_m(t)$:
\begin{align}
    y_m(t) = Q(p(t-D)) \label{eq:medicao_atrasada}.
\end{align}

\begin{assum}
The maximizer $p^* = [x^*, y^*]^\top \in \mathbb{R}^2$ and the maximum value $Q^* \in \mathbb{R}$ of the nonlinear map defined in \eqref{eq:mapa_Q} are static and unknown parameters.
\end{assum}

\begin{assum}
The measurement delay $D > 0$ is constant and known.
\end{assum}

\begin{assum}
The system operates without global positioning coordinates. The only variables available to the ESC algorithm and the predictor mechanism are the measured signal $y_m(t)$ defined in \eqref{eq:medicao_atrasada} and the time history of the control signals $u(\tau) = [u_x(\tau), u_y(\tau)]^\top  \in \mathbb{R}^2$ applied over the interval $\tau \in [t-D, t]$.
\end{assum}

\noindent \textit{Probing and Demodulation Signals for Averaging-Based Estimation}

Define $\hat{p}(t) = [\hat{x}_p(t), \hat{y}_p(t)]^\top \in \mathbb{R}^2$ as the current estimate of the source position. The position estimation error is given by
\begin{align}
    \tilde{p}(t) &= \hat{p}(t) - p^* \label{eq:erro_estimativa}.
\end{align}

To perform the search in the $2\text{D}$ plane, high-frequency dither signals $S(t) = [S_x(t), S_y(t)]^\top \in \mathbb{R}^2$ are injected. To compensate for the measurement latency $D$, the dither signals are introduced with a phase advance $+D$:
\begin{align}
    S(t) &= \begin{bmatrix} 
        a_1 \sin(\omega_1 (t + D)) \\[4pt] 
        a_2 \cos(\omega_2 (t + D)) 
    \end{bmatrix} \label{eq:vetor_dither},
\end{align}
where $a_1, a_2 > 0$ are the amplitudes and $\omega_1 \neq \omega_2$ are distinct and incommensurate dither frequencies.

The actual position of the point $p(t)$ at the current time is the combination of the estimate and the dither signal, i.e., $p(t) = \hat{p}(t) + S(t)$. Consequently, the delayed point perceived through the measurement $y_m(t) = Q(p(t-D))$ is given by $p(t-D) = \hat{p}(t-D) + S(t-D)$.

\newpage
The estimates of the local gradient $\hat{G}(t) \in \mathbb{R}^2$ and Hessian matrix $\hat{H}(t) \in \mathbb{R}^{2 \times 2}$ are obtained by demodulating the measured signal $y_m(t)$:
\begin{align}
    \hat{G}(t) &= M(t) y_m(t) \label{eq:gradiente_demodulado}, \\[4pt]
    \hat{H}(t) &= N(t) y_m(t) \label{eq:estimativa_hessiana},
\end{align}
where the gradient demodulation vector $M(t)$ and the Hessian demodulation matrix $N(t)$ are defined by
\begin{align}
    M(t) &= \begin{bmatrix} 
        \frac{2}{a_1} \sin(\omega_1 t) \\[4pt] 
        \frac{2}{a_2} \cos(\omega_2 t) 
    \end{bmatrix} \label{eq:vetor_demodulacao}, \\[6pt]
    N(t) &= \begin{bmatrix} 
        -\frac{8}{a_1^2}\cos(2\omega_1 t) & 0 \\[6pt] 
        0 & -\frac{8}{a_2^2}\cos(2\omega_2 t) 
    \end{bmatrix} \label{eq:matriz_demodulacao_2a_ordem}.
\end{align}

By applying averaging theory over the common fundamental period $T = \frac{2\pi}{\gcd(\omega_1, \omega_2)}$, moving integration of the signals yields the averaged components of the system:
\begin{align}
    \hat{H}_{av} &= \frac{1}{T} \int_{t-T}^{t} N(\sigma) y_m(\sigma) \, d\sigma = H, \qquad H:= \begin{bmatrix} -q_{x} & 0 \\ 0 & -q_{y} \end{bmatrix}<0  \label{eq:media_hessiana_integral}, \\[4pt]
    \hat{G}_{av}(t) &= \frac{1}{T} \int_{t-T}^{t} M(\sigma) y_m(\sigma) \, d\sigma = H \tilde{p}_{av}(t-D) \label{eq:media_gradiente_integral},
\end{align}
where $H \in \mathbb{R}^{2 \times 2}$ is the Hessian matrix of the quadratic map \eqref{eq:mapa_Q}. Equation \eqref{eq:media_gradiente_integral} confirms that the average of the demodulated gradient at the current time $t$ provides the past position error $H \tilde{p}_{av}(t-D)$.\\

\noindent \textit{Predictor for Delay Compensation.}

From the averaged model, the time dynamics of the delayed position error and the averaged gradient are given by
\begin{align}
    \frac{d\tilde{p}_{av}(t - D)}{dt} &= u_{av}(t - D) \label{eq:dinamica_erro_medio}, \\
    \frac{d\hat{G}_{av}(t)}{dt} &= H u_{av}(t - D) \label{eq:dinamica_gradiente_medio},
\end{align}
where $u_{av}(t) \in \mathbb{R}^2$ represents the vector of virtual averaged control signals for the point $p(t)$.

To compensate for the measurement delay $D$, the future state $\hat{G}_{av}(t+D)$ is predicted using the variation-of-constants formula applied to \eqref{eq:dinamica_gradiente_medio}:
\begin{align}
    \hat{G}_{av}(t + D) &= \hat{G}_{av}(t) + H \int_{t-D}^{t} u_{av}(\sigma) \, d\sigma \label{eq:predicao_estado_futuro}.
\end{align}

The ideal predictive control law in the averaged domain is specified as $u_{av}(t) = k \hat{G}_{av}(t + D)$, where $k > 0$ is the control gain. Substituting \eqref{eq:predicao_estado_futuro} yields the implementable feedback law:
\begin{align}
    u_{av}(t) &= k \left[ \hat{G}_{av}(t) + H \int_{t-D}^{t} u_{av}(\sigma) \, d\sigma \right] \label{eq:lei_controle_preditiva_media}.
\end{align}

Applying \eqref{eq:lei_controle_preditiva_media} to \eqref{eq:dinamica_erro_medio}, the closed-loop averaged error dynamics become
\begin{align}
    \frac{d\tilde{p}_{av}(t)}{dt} &= k H \tilde{p}_{av}(t), \quad \forall t \ge D \label{eq:sistema_fechado_medio},
\end{align}
which are exponentially stable at the origin $\tilde{p}_{av} = 0$, since $H < 0$.

To enable implementation in the continuous-time system and approximate the averaged signals $\hat{G}_{av}(t)$ and $\hat{H}_{av}(t)$, a low-pass filter with cutoff frequency $c > 0$ is applied. The control signal generated by the predictor, $u(t) \in \mathbb{R}^2$, is governed by the differential equation:
\begin{align}
    \dot{u}(t) &= -c u(t) + c k \left[ \hat{G}(t) + \hat{H}(t) \int_{t-D}^{t} u(\tau) \, d\tau \right] \label{eq:preditor_filtrado_tempo}.
\end{align}

Finally, the Cartesian control vector $u(t) = [u_x(t), u_y(t)]^\top$ generated by the predictor \eqref{eq:preditor_filtrado_tempo} is converted into the actual unicycle velocities $v(t)$ and $\omega(t)$ by inverting the decoupling matrix \eqref{eq:matriz_decouplage}:
\begin{align}
    \begin{bmatrix} v(t) \\[4pt] \omega(t) \end{bmatrix} &= T^{-1}(\theta(t)) u(t) = \begin{bmatrix} \cos\theta(t) & \sin\theta(t) \\[4pt] -\frac{1}{l}\sin\theta(t) & \frac{1}{l}\cos\theta(t) \end{bmatrix} \begin{bmatrix} u_x(t) \\[4pt] u_y(t) \end{bmatrix} \label{eq:mapeamento_atuadores}.
\end{align}

\begin{figure}
	\centering
	\includegraphics[width=12cm]{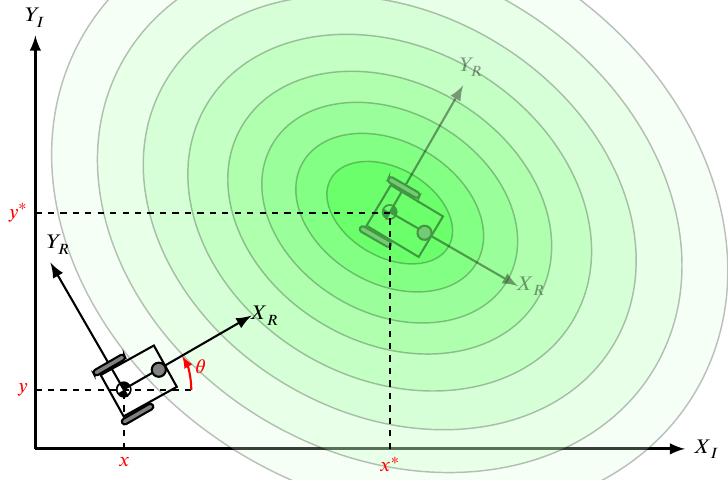}
	\caption{Source-seeking geometry for a nonholonomic mobile robot. The robot pose is characterized by its global position $(x,y)$ and heading angle $\theta$ between the robot-fixed axis $X_R$ and the global axis $X_I$. The unknown source is located at $(x^*,y^*)$, with the shades of green depicting the spatial gradient of the source signal field.}
	\label{sourceseeking}
\end{figure}

\newpage
\subsection{Further Developments on Extremum Seeking} 


\noindent \textit{Extremum Seeking for Distinct Families of PDEs.}

\begin{figure}
	\centering
	\includegraphics[width=10cm]{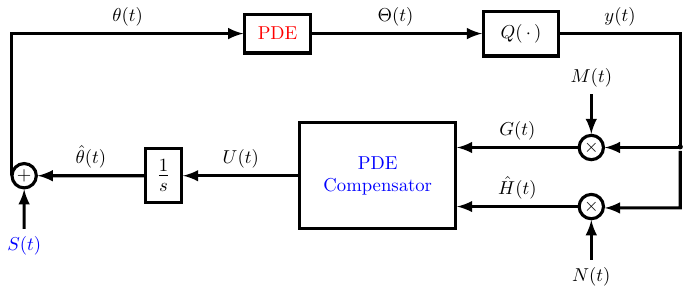}
	\caption{A general block diagram for implementation of ESC design for PDE connections (\textcolor{red}{in red}) at the input of nonlinear convex maps $Q(\cdot)$.
	Although the multiplicative perturbation signals $M(t)$ and $N(t)$ are the same as in the classical ESC designs \cite{KrsticWang2000,GKN:2012}, the additive dither $S(t)$ (\textcolor{blue}{in blue}) must be redesigned using the trajectory generation paradigm \cite[Ch.~12]{krstic2008boundary} and the application of an adequate boundary control law (\textcolor{blue}{in blue}) for PDE compensation is necessary.}
	\label{ch12.fig:cascade}
\end{figure}

In order to show that the proposed ESC approach for infinite-dimensional systems is general and applicable to a wider class of PDEs,
we have formulated in Table~\ref{tabelademerda} the stabilizing boundary control (BC) law $U(t)$ and given the explicit solutions to the trajectory generation problem of $S(t)$  
for five other classes of distributed parameter systems \cite{krstic2008boundary}: (a) reaction-advection-diffusion (RAD)
PDEs \cite{Alcos:2019,IJACSP2020}, (b) wave equations \cite{OK:2019,ASME2021},
(c) hyperbolic transport PDEs ---for constant delays \cite{Tiago},
(d) time-varying delays \cite{ORDK:2018}, and (e) distributed delays \cite{TOK:2020,AUT2023}.


\begin{table*}[ht]
	\caption{ESC for Distinct Classes of PDE Systems.}
	\centering
	\resizebox{15cm}{!}{
	\begin{tabular}{||c||l||}
		\hline
  \hline
		\\{\LARGE\textbf{RAD Equation}}    &\\                     
		              &
									{\LARGE$\textcolor{blue}{\text{PDE}:} ~\partial_t\alpha(x,t)\!=\!\epsilon\partial_{xx}\alpha(x,t)\!+\!
									b\partial_x\alpha(x,t)\!+\! \lambda\alpha(x,t)\,, \quad x\in[0,1]$} \\
		              &{\LARGE$\textcolor{blue}{\text{BC (Dirichlet):}} ~U(t)\!=\! \frac{c}{s+c}\left\lbrace k
									e^{-\frac{b}{2\epsilon}}\left[ \gamma(1)G(t)\!+\!\hat{H}(t)\int_{0}^{1}e^{\frac{b}{2\epsilon}\sigma} m(1\!-\!\sigma)
									u(\sigma,t)d\sigma \right]\right\rbrace,$}\\
									 &{\LARGE$\gamma(x) = \cosh\left(\sqrt{\frac{\xi}{\epsilon}}x\right)+\frac{b}{2\epsilon}\sqrt{\frac{\epsilon}{\xi}}									\sinh\left(\sqrt{\frac{\xi}{\epsilon}}x\right), \quad \xi:=b^2/(4\epsilon)-\lambda\ge 0\,, k>0$}\\
									&{\LARGE$m(x-\sigma) = \frac{1}{\epsilon}\sqrt{\frac{\epsilon}{\xi}}\sinh\left(\sqrt{\frac{\xi}{\epsilon}}
									(x-\sigma)\right), \quad \epsilon>0\,, b\geq 0\,, \lambda \geq 0$}\\
								  &{\LARGE$\textcolor{blue}{\text{Trajectory Generation}:} ~S(t) = e^{-\frac{b}{2\epsilon}}\sum_{k=0}^{\infty}
									\frac{a_{2k}(t)}{(2k)!}+\frac{b}{2\epsilon}\frac{a_{2k}(t)}{(2k+1)!},$}\\
									&{\LARGE$a_{2k} := \frac{a}{\epsilon^k}\sin(\omega t) \sum_{n=0}^{k}\binom{k}{2n}\xi^{k-2n}\omega^{2n}
									+ \frac{a}{\epsilon^k}\cos(\omega t)\sum_{n=0}^{k}\binom{k}{2n+1}\xi^{k-2n-1} \omega^{2n+1}$}\\
		&\\
		\hline
  \hline
		\\{\LARGE\textbf{Wave Dynamics}} &\\
									                                 &
									{\LARGE$\textcolor{blue}{\text{PDE}:} ~\partial_{tt}\alpha(x,t)=\partial_{xx}\alpha(x,t), \quad x\in[0,D]$} \\
		              &{\LARGE$\textcolor{blue}{\text{BC (Neumann):}} ~U(t)=\frac{c}{s+c}
									\left\{c\left[k \hat{H}(t)u(D,t) - \partial_t u(D,t)\right] + \rho(D) G(t) + 
									\phantom{\int_{0}^{D}} \right.$} \\
									&{\LARGE$\left.  \hat{H}(t)\int_{0}^{D} \rho(D-\sigma) \partial_t u(\sigma,t) d\sigma\right\}$, $\quad 
									\rho(s)=k[0 \ \ I]e^{As}[0 \ \ I]^T, \quad A=\begin{pmatrix}
		                                    0 & 0 \\ I & 0 \\
	                                      \end{pmatrix}$}\\
								  &{\LARGE$\textcolor{blue}{\text{Trajectory Generation}:} ~S(t) = a
									\cos(\omega D) \sin(\omega t)$}\\
		&\\
		\hline
  \hline 
		\\{\LARGE\textbf{Constant Delay}}     &\\       &
									{\LARGE$\textcolor{blue}{\text{PDE}:} ~\partial_{t}\alpha(x,t)=\partial_{x}\alpha(x,t), \quad x\in[0,D]$} \\
		              &{\LARGE$\textcolor{blue}{\text{BC (Dirichlet)}:} ~\textcolor{black}{U(t)=\frac{c}{s+c}\left\{k 
									\left[G(t)+\hat{H}(t)\int_{0}^{D} u(\sigma,t) d\sigma \right]\right\}}$}\\
								  &{\LARGE$\textcolor{blue}{\text{Trajectory Generation}:} ~\textcolor{black}{S(t) = a\sin(\omega (t+D))}$}\\
		&\\
		\hline
  \hline
		\\{\LARGE\textbf{Variable Delay}}      &\\       &
									{\LARGE$\textcolor{blue}{\text{PDE}:} ~\partial_{t}\alpha(x,t)= \pi(x,t)\partial_{x}\alpha(x,t),
									\quad x\in[0,1], \quad \pi(x,t)=\frac{1+x[\frac{d(\phi^{-1}(t))}{dt}-1]}{\phi^{-1}(t)-t}$} \\
		              &{\LARGE$\textcolor{blue}{\text{BC (Dirichlet):}} ~\textcolor{black}{U(t)=\frac{c}{s+c}\left\{k 
									\left[G(t)+\hat{H}(t)\int_{0}^{1} u(\sigma,t)\left(\phi^{-1}(t)-t\right) d\sigma \right]\right\}}$}\\
								  &{\LARGE$\textcolor{blue}{\text{Trajectory Generation}:} 
									~\textcolor{black}{S(t) = a\sin(\omega t)},\quad \phi(t):= t-D(t)$}\\ 
									&{\LARGE$\textcolor{blue}{\text{Demodulation}:} ~M(t) = \frac{2}{a}\sin(\omega (t-D(t)))\,, ~N(t) = -\frac{8}{a^2}\cos(2\omega (t-D(t)))$}\\
									
		&\\
		\hline
  \hline
		\\{\LARGE\textbf{Distributed Delay}}       &\\              &
									{\LARGE$\textcolor{blue}{\text{PDE}:}  ~\partial_{t}\alpha(x,t)=\partial_{x}\alpha(x,t), \quad x\in[0,D], 
									\quad y=Q\left( \int_{0}^{D} \Theta(t-\sigma) d\beta(\sigma)  \right)$} \\
		              &{\LARGE$\textcolor{blue}{\text{BC (Dirichlet):}} ~\textcolor{black}{U(t)=\frac{c}{s+c}\left\{k 
									\left[G(t)+\hat{H}(t)\int_{0}^{D} (1-\beta(\sigma))u(D-\sigma,t) d\sigma \right]\right\}}$}\\
								  &{\LARGE$\textcolor{blue}{\text{Trajectory Generation}:} 
									~\textcolor{black}{S(t) = \frac{a}{\gamma(\omega)}\int_{0}^{D}\sin(\omega (t+\xi))d\beta(\xi)}$}\\							
		&\\
		\hline
  \hline
	\end{tabular}
	}
\label{tabelademerda}
\end{table*}


Regarding the \textcolor{black}{transport and diffusion cases}, the term $u(x,t)$ which appears in $U(t)$ of Table~\ref{tabelademerda}, is the state of the infinite-dimensional system corresponding to a copy of the PDE model of actuator dynamics. 

In \cite{PDE_cascades_SCL}, differently from what has been done that has dealt with PDEs at the input into an unknown map, and in which we have already advanced from transport PDEs to reaction-advection-diffusion PDEs to wave PDEs according to Table~~\ref{tabelademerda}, the authors consider one last configuration in which the input pathway to the map contains a cascade of PDEs from distinct classes. 
There, we deal with PDEs with input delays, such as, the notorious problem of a wave PDE with an input delay where, if the delay is left uncompensated, an arbitrarily short delay destroys the closed-loop stability. Then, we move forward to an even more challenging class of problems for parabolic-hyperbolic cascades of PDEs, focusing on a heat equation at the input of a wave PDE.  
The treatment of such systems with PDE-PDE cascades is again performed by means of boundary control. 
Local exponential stability and convergence to a small neighborhood of the unknown extremum point are  guaranteed by using a backstepping transformation and averaging in infinite dimensions. 

PDE-PDE cascades have a great deal in common with PDE-ODE cascades. For instance, a cascade of a delay into a PDE is a much more generalized version of an integrator with an input delay. 
However, while in a delay-integrator cascade the design can be pursued within the predictor feedback framework, with backstepping just employed for an interpretation and for analysis, in delay-PDE cascades a predictor for a PDE is too complicated of a mathematical object to be of value. Instead, the design is pursued entirely by the backstepping approach. Similarly, while the heat-integrator cascade became familiar, 
and was dealt with through the backstepping design, a more general backstepping design is applied to a heat-wave PDE cascade in a part of the reference \cite{PDE_cascades_SCL}.

\begin{figure}
\centering
\includegraphics[width=12cm]{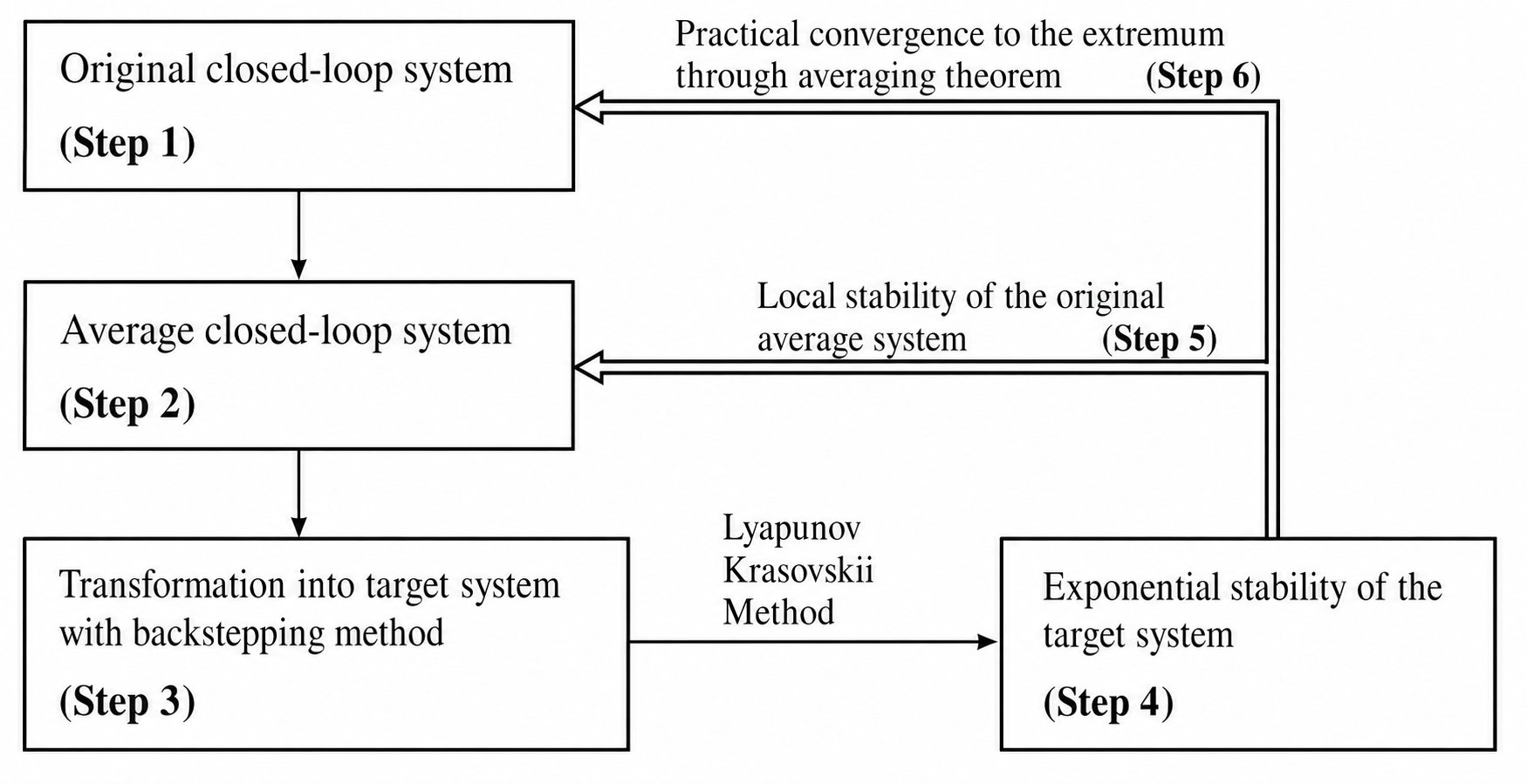}
\caption{\textcolor{black}{Structure of the stability proof for the closed-loop system.}} \label{ch9.fig:structure_proof}
\end{figure}

Although our goal in this paper is to avoid numbing the reader with lengthy proofs, it is worth providing a general picture on how we carry out the steps to prove the stability results of the ESC feedback loop in the presence of PDEs in Figure~\ref{ch12.fig:cascade}. Figure~\ref{ch9.fig:structure_proof} shows the structure of the proof, divided into six main steps. We take advantage of this opportunity to highlight that our analysis presents a carefully constructed sequence of analytical steps, a predictor-based infinite-dimensional backstepping transformation, a synthesis of a Lyapunov functional (rather than small-gain analysis), and computation of a Lyapunov estimate, for the overall infinite-dimensional system with nonlinearities, stochastic perturbations, and distributed delays. The analysis process involving so many steps has a large number of possible permutations---all of which, \textit{but one} would be wrong. We show how to properly sequence the steps of averaging, backstepping, and Lyapunov functional analysis to prove stability. This ``analysis pathway'' will serve the needs of future researchers who deal with stochastic extremum seeking under delays. The complete details can be found in the author's book \cite{OliveiraKrstic2022}. 
%

For the distinct cases involving ESC plus PDEs in Table~~\ref{tabelademerda}, we can list some selected applications:

\textbf{1)} \textit{Traffic Control} for linearized Lighthill-Whitham-Richards (LWR) macroscopic PDE models transformed into constant delays \cite{YKOK:2021,Tiago}; 
\textbf{2)} \textit{Optimal Oil Drilling Control} with ESC for wave models \cite{AK:2019}; 
\textbf{3)} \textit{Deep-Sea Cable-Actuated Source Seeking} modeled by wave PDEs with Kelvin-Voigt damping \cite{L-CSS:2023}; 
\textbf{4)} \textit{Additive manufacturing} modeled by the \textit{Stefan PDE} \cite{KDK:2019,TDS2022}; 
\textbf{5)} \textit{Bioreactors} considering ESC for models described by parabolic PDEs (reaction-diffusion equations), see \cite{DOCHAIN:2006,DOCHAIN:2008};  
\textbf{6)} \textit{Light-Source Seeking} with infinite-dimensional models represented by Euler-Bernoulli beams equations \cite[Ch.~8]{krstic2008boundary}; and  
\textbf{7)} \textit{Neuromuscular Electrical Stimulation} (NMES) problem for ESC with time-varying delays \cite{PAZ:2019,ORDK:2018}. \\

\noindent \textit{Event-Triggered Extremum Seeking Control.} 

\begin{figure}
\centering
\includegraphics[width=15cm]{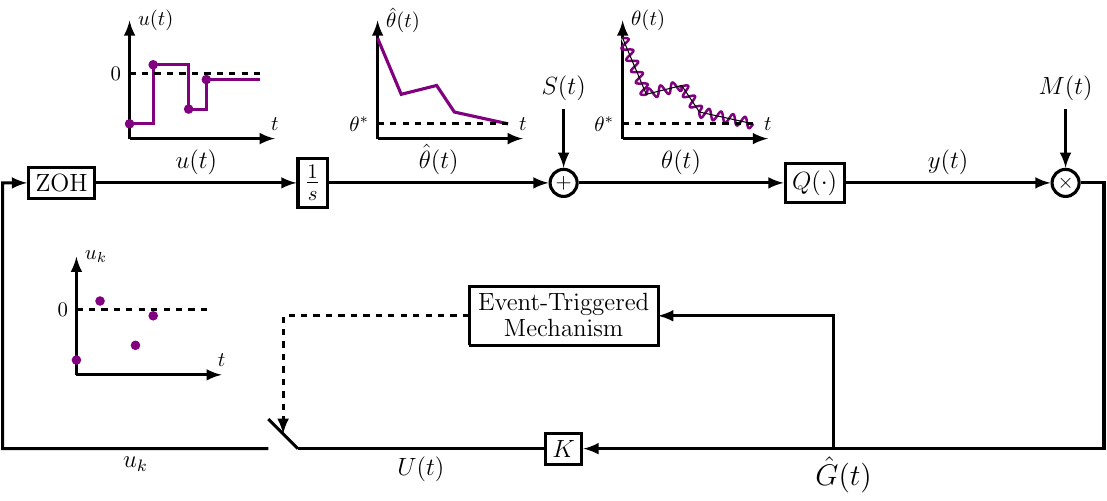}
\caption{\textcolor{black}{Event-triggered extremum-seeking architecture, illustrating the transition from continuously updated ESC toward resource-aware optimization, in which the gradient-based control action is updated only when needed. Adapted from \cite{RodriguesEtAl2025}.}} \label{ET-ESC}
\end{figure}

As extremum seeking has moved from its classical continuous-time formulation toward digitally implemented and networked systems, a natural question has emerged: \textit{does the optimization loop really need to be updated continuously?}
Event-triggered extremum seeking control (ET-ESC) provides a negative answer by bringing the resource-aware philosophy of event-triggered control into the classical perturbation-based ESC architecture \cite{RodriguesEtAl2025}. As illustrated in Figure~\ref{ET-ESC}, the usual dither and demodulation machinery is retained, but the gradient-feedback control action is updated only when an event-triggering condition indicates that a new update is needed. The resulting control signal is piecewise constant and, because the ESC dynamics contains an integrator, the parameter estimate becomes piecewise linear between events. In this way, optimization performance can be traded against communication, computation, and actuation effort rather than being tied to a prescribed periodic sampling rate. Importantly, this event-driven viewpoint removes the upper-bound restriction on inter-sampling times encountered in conventional sampled-data ESC.

Static and dynamic event-triggering mechanisms were developed and analyzed through Lyapunov arguments combined with averaging theory for discontinuous systems, establishing exponential stability of the averaged dynamics, practical convergence of the actual trajectories to the unknown extremum, and exclusion of Zeno behavior through a positive minimum inter-event time. A periodic event-triggered version (PET-ESC) further moves this philosophy toward digital implementation by evaluating the triggering condition only at periodic instants while retaining aperiodic control updates. Viewed from a broader perspective, these developments introduce \textbf{resource awareness as a new design dimension in extremum seeking}: the question is no longer only how rapidly and accurately an unknown optimum can be reached, but also how much sensing, communication, computation, and actuation must be spent to reach it. This viewpoint points naturally toward future ESC architectures for networked, embedded, distributed, and large-scale optimization systems, where the intelligent allocation of resources may become as important as the optimization objective itself. 

\newpage
\noindent \textit{Extremum Seeking under Saturation Constraints.} 

\begin{figure}
\centering
\includegraphics[width=10cm]{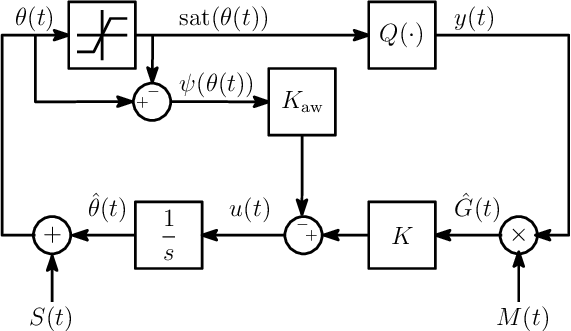}
\caption{\textcolor{black}{Extremum-seeking architecture under actuator saturation. An anti-windup feedback loop accounts for the difference between the commanded and admissible map inputs, illustrating saturation as a physical limitation to be compensated. Adapted from \cite{SilvaEtAl2026}.}} \label{SAT1}
\end{figure}

As extremum seeking moves from idealized optimization loops toward physical systems, another fundamental question arises: \textit{what happens when the optimizer asks for more than the system can deliver?} Actuators, commands, and adaptation rates are inevitably bounded in practice, and ignoring such limitations may degrade performance or even compromise closed-loop stability. Saturated extremum seeking addresses this gap by explicitly incorporating hard constraints into the ESC feedback architecture \cite{SilvaEtAl2026}. Two complementary viewpoints emerge. In the first, saturation represents a \textbf{physical limitation to be compensated}: when the input commanded by the optimizer exceeds the admissible range, an anti-windup mechanism uses the mismatch between the commanded and saturated inputs to modify the adaptation dynamics, as illustrated in Figure~\ref{SAT1}. Unlike earlier approaches largely restricted to scalar problems or decentralized anti-windup compensation, the multivariable formulation permits the joint design of full, potentially non-diagonal feedback and anti-windup gains while accounting for uncertainty in the unknown Hessian.

\begin{figure}
\centering
\includegraphics[width=10cm]{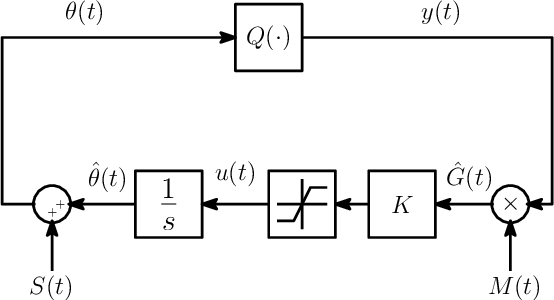}
\caption{\textcolor{black}{Extremum-seeking architecture under gradient saturation. By saturating the gradient-feedback signal before integration, the update rate of the parameter estimate is explicitly bounded, illustrating saturation as a design mechanism rather than merely a physical limitation. Adapted from \cite{SilvaEtAl2026}.}} \label{SAT2}
\end{figure}

A complementary viewpoint arises when saturation is deliberately placed in the gradient-feedback path, as shown in Figure~\ref{SAT2}. Here, saturation is no longer merely a limitation to overcome but becomes a \textbf{design mechanism to be exploited}: by bounding the signal entering the ESC integrator, the update rate of the parameter estimate is explicitly bounded. This provides a formulation of classical perturbation-based ESC with bounded update rates while preserving convergence toward a neighborhood of the unknown extremum. For both architectures, sector representations, polytopic descriptions of the uncertain Hessian, Lyapunov arguments, and averaging theory for systems with Lipschitz right-hand sides provide constructive stability guarantees and controller-design conditions. More broadly, these two interpretations bring \textbf{constraint awareness into extremum seeking}: optimization can no longer be concerned only with where the optimum lies and how rapidly it is reached, but must also respect how far and how fast the physical system is allowed to move. This viewpoint points toward ESC architectures in which actuator limits, safety envelopes, rate constraints, and other physical restrictions become intrinsic elements of the optimization design rather than implementation details considered afterward. \\


\newpage
\noindent \textit{Resilient Extremum Seeking under Cyber-Physical Attacks.} 

\begin{figure}
\centering
\includegraphics[width=8.9cm]{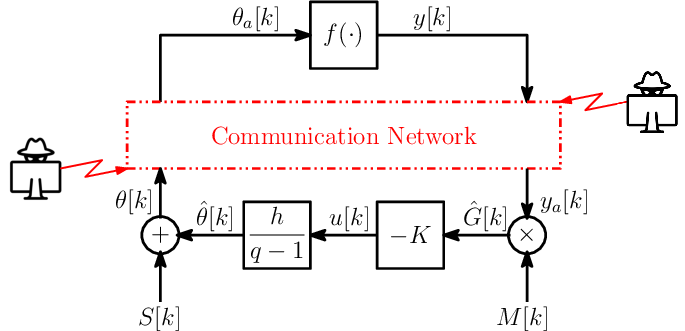}
\caption{\textcolor{black}{Networked extremum-seeking architecture under denial-of-service attacks. Communication between the optimizer and the physical process may be intermittently disrupted in both the measurement and control channels, motivating resilient ESC strategies that preserve convergence despite unavailable data. In the corresponding discrete-time configuration, we consider that any sequence $f=\{f[k]\}_{k\in\mathbb{Z}}$, $f[k] \in \mathbb{R}^{n}$, is obtained by sampling the original continuous function $f_a(t)$ with period $h>0$, {\it i.e.}, $f[k] = f_a(kh)$. The {\bf forward} and {\bf backward shift operators} are $qf[k] = f[k+1]$ and $q^{-1}f[k] = f[k-1]$; in the $z$-domain, $q$ is identified with the complex variable $z \in \mathbb{C}$. Adapted from \cite{TeixeiraEtAl2026}.}} \label{DOS}
\end{figure}

As extremum seeking becomes increasingly embedded in networked and cyber-physical systems, a new question arises: \textit{can the optimizer continue to operate when the communication infrastructure itself is under attack?} This issue motivates resilient networked (discrete-time) extremum seeking against denial-of-service (DoS) attacks \cite{TeixeiraEtAl2026}. As illustrated in Figure~\ref{DOS}, the feedback loop is closed through a communication network: measurements of the performance output must reach the ESC algorithm, while the optimized input must be transmitted back to the physical process. A DoS attacker can interrupt either exchange, temporarily depriving the optimizer of fresh measurements or preventing newly computed inputs from reaching the process. A mixed mitigation strategy is adopted in which the unavailable output is replaced by zero, whereas the previously received input is held during communication outages. Consequently, the ESC dynamics alternates between a \textbf{stable optimization mode}, when communication succeeds, and a \textbf{neutral mode}, when attacks block packet transmissions and the optimization error is essentially frozen rather than amplified.

\newpage
Resilience is established for two complementary descriptions of cyber attacks. For deterministic DoS, the attacker may interrupt communication over arbitrary intervals provided that the overall fraction of time under attack remains bounded; practical exponential convergence is guaranteed whenever successful communication persists for a nonvanishing fraction of time. For probabilistic DoS, packet failures are modeled stochastically through their attack-success probability, leading to almost-sure and probabilistic stability guarantees as long as communication is not permanently suppressed. Viewed more broadly, these results introduce \textbf{cyber resilience as another design dimension in extremum seeking}: future optimization loops must not only determine where the optimum lies, use resources judiciously, and respect physical constraints, but also tolerate interruptions and adversarial manipulation of the information required to reach it. This perspective becomes increasingly important as ESC migrates toward interconnected autonomous systems, distributed optimization platforms, industrial networks, and other cyber-physical environments in which communication can no longer be assumed continuously available or trustworthy. \\

\noindent \textit{Sliding-Mode-Based Extremum Seeking Control.} 

\begin{figure}
\centering
\includegraphics[width=7cm]{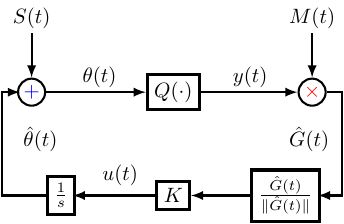}
\caption{\textcolor{black}{Gradient-based unit-vector extremum-seeking architecture. Normalization of the estimated-gradient feedback replaces the conventional proportional action by a relay-type unit-vector control law, yielding finite-time convergence of the averaged optimization dynamics. Adapted from \cite{LuoEtAl2026}.}} \label{Gradient_UVC}
\end{figure}

As extremum seeking continues to evolve toward faster optimization schemes, a fundamental question arises: \textit{must convergence to the extremum remain asymptotic, or can the optimization dynamics reach its target in finite time?} Sliding-mode-based extremum seeking provides a natural route toward the latter by importing discontinuous control mechanisms into the classical perturbation-based ESC architecture \cite{LuoEtAl2026}. In the gradient-based unit-vector ESC (UV-ESC) depicted in Figure~\ref{Gradient_UVC}, the conventional proportional action on the estimated gradient is replaced by a normalized, relay-type feedback law. The resulting control action retains the direction of the estimated gradient while normalizing its magnitude, producing a discontinuous feedback at the target. This seemingly simple modification fundamentally changes the convergence mechanism: whereas conventional gradient-based ESC with proportional feedback yields exponential stability of the averaged dynamics, the unit-vector formulation drives the averaged gradient estimate to zero in finite time, while the actual ESC trajectories converge to an arbitrarily small neighborhood of the unknown optimizer for sufficiently small dither amplitudes and sufficiently high probing frequencies.

\begin{figure}
\centering
\includegraphics[width=14cm]{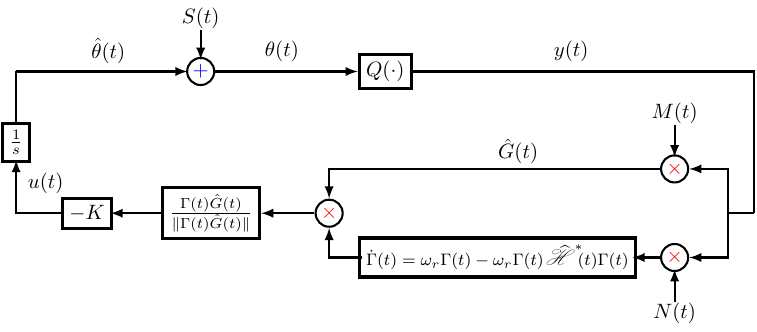}
\caption{\textcolor{black}{Newton-based unit-vector extremum-seeking architecture. A Riccati-based inverse-Hessian estimator generates a Newton direction whose normalized feedback combines finite-time convergence with a convergence rate independent of the unknown Hessian, with the averaging-based estimate ~~~~$\widehat{\!\!\!\!\!\!\!\!\mbox{\calligra H}}^{~~\ast}(t)=N(t)y(t)$. Adapted from \cite{LuoEtAl2026}.}} \label{Newton_UVC}
\end{figure}

The same philosophy extends naturally to Newton-based ESC, as illustrated in Figure~\ref{Newton_UVC}. Here, second-order information is estimated through demodulation and a Riccati-based estimator of the inverse Hessian, while unit-vector feedback is applied to the resulting Newton direction. This combination retains finite-time convergence while overcoming the dependence of the gradient-based convergence rate on the unknown Hessian, thereby accelerating the transient response. More broadly, these developments introduce \textbf{finite-time convergence as another design dimension in extremum seeking}, connecting model-free real-time optimization with the rich machinery of variable-structure and sliding-mode control. First-order unit-vector feedback represents only one entry point into this largely unexplored interface: higher-order, adaptive, and integral sliding modes, as well as fixed- and prescribed-time mechanisms, suggest a broader family of ESC algorithms in which not only the destination of the optimization process, but also \textbf{how—and how quickly—it is reached}, becomes an explicit part of the controller design. \\


Taken together, these developments illustrate a broader transformation in the objectives of extremum-seeking design. Future ESC algorithms will be expected not only to locate an unknown optimum, but also to determine \textbf{how fast it should be reached, how many resources should be spent in reaching it, which physical constraints must be respected along the way, and how optimization can be sustained under adverse or even malicious operating conditions}. Event-triggered, saturated, cyber-resilient, and sliding-mode-based ESC provide complementary steps in this direction, bringing resource awareness, constraint awareness, resilience, and finite-time convergence into the optimization loop. These advances suggest a transition from classical model-free optimization toward a new generation of \textit{resource-, constraint-, resilience-, and convergence-aware extremum-seeking architectures}, in which performance is measured not only by proximity to the optimum, but also by the efficiency, feasibility, robustness, and speed with which that optimum is pursued.

\newpage
\section{Where the Field Is Going Next}

\noindent \textit{Industry 4.0 and the Age of Autonomous Optimization.} 

The emergence of Industry 4.0 has only reinforced the relevance of extremum seeking. Modern industrial systems increasingly rely on networks of sensors, actuators, communication technologies, and embedded intelligence operating continuously under changing conditions. In such environments, accurate mathematical models quickly become outdated, while operating conditions may vary from one moment to the next.

Extremum seeking offers a natural solution to this challenge. Instead of relying exclusively on detailed system models, it continuously learns from measured performance, allowing industrial processes to adapt autonomously to changing operating conditions while maintaining mathematical guarantees of stability and convergence.

In this sense, extremum seeking represents an important bridge between classical feedback control and the emerging paradigm of intelligent autonomous systems. It combines the reliability of model-based control theory with the adaptability required by modern cyber-physical systems, making it particularly attractive for the next generation of industrial technologies. \\

\noindent \textit{One Principle, Countless Applications.} 

One of the most fascinating aspects of extremum seeking is that its underlying principle is almost entirely independent of the physical system to which it is applied. Whether the objective is to maximize power generation, minimize fuel consumption, improve communication quality, or optimize a biomedical treatment, the optimization algorithm itself remains essentially unchanged. Only the dynamical system surrounding it is different.

This universality explains why extremum seeking has found applications across an extraordinary range of disciplines. In modern automobiles, it can continuously adjust engine operating conditions to improve fuel efficiency and reduce emissions. Wind turbines use similar principles to maximize power extraction under continuously changing wind conditions. Industrial processes employ extremum seeking to optimize temperature, pressure, flow rates, and chemical reactions without requiring precise process models.

The same methodology has also proved remarkably effective in areas that might initially appear far removed from classical control engineering. Examples include telecommunications, biological reactors, additive manufacturing, autonomous marine vehicles, robotics, renewable energy systems, smart buildings, and neuromuscular electrical stimulation. In every case, the controller performs the same fundamental task: it continuously learns, through interaction with the physical system, how to improve performance despite uncertainty.

This remarkable transferability is uncommon in engineering. Many control algorithms are designed for highly specific classes of systems. Extremum seeking, by contrast, has demonstrated an unusual ability to migrate across disciplines while preserving both its conceptual simplicity and its mathematical foundations.

In retrospect, this universality may be one of the principal reasons for its continued expansion over the past two decades. \\

\noindent \textit{Beyond Optimization: A New View of Learning.} 

Historically, extremum seeking has been described as a real-time optimization algorithm. While this description is certainly correct, it no longer captures the full scope of what the field has become.

At its core, extremum seeking is a mechanism for \textbf{learning through feedback}. The controller deliberately explores its environment, extracts information from the measured response, and immediately incorporates that information into future decisions. Learning, optimization, and control are therefore not separate processes but different manifestations of the same closed-loop dynamics.

Seen from this perspective, extremum seeking belongs to a much broader scientific movement that seeks to endow engineered systems with increasing levels of autonomy. Modern cyber-physical systems are expected not only to regulate their behavior, but also to adapt, optimize, cooperate, and continue improving while operating in uncertain and evolving environments.

This broader interpretation helps explain why extremum seeking has recently begun interacting with fields such as distributed optimization, game theory, reinforcement learning, machine learning, and autonomous decision-making. Although these communities often employ different terminology and mathematical tools, they ultimately pursue the same fundamental objective: enabling dynamical systems to improve their own performance through experience. \\

\newpage
\noindent \textit{Extremum Seeking and Artificial Intelligence: Partners Rather Than Competitors.} 

The rapid progress of artificial intelligence has naturally prompted questions about the future role of classical feedback optimization techniques. At first glance, one might imagine that reinforcement learning or other data-driven methods will eventually replace algorithms such as extremum seeking. We believe the opposite is more likely.

Artificial intelligence and extremum seeking address complementary aspects of autonomous decision-making. Machine learning excels at discovering complex patterns from large datasets and constructing predictive models of uncertain environments. Extremum seeking, on the other hand, provides rigorous mechanisms for stable real-time optimization directly within the feedback loop, often requiring only current performance measurements rather than extensive offline training.

Rather than competing, these paradigms are likely to become increasingly integrated. Machine learning may provide richer representations of the environment, while extremum seeking supplies mathematically grounded adaptation mechanisms that preserve robustness and closed-loop stability. The result could be intelligent control architectures capable of combining data-driven prediction with provable real-time optimization.

This convergence represents one of the most promising directions for the next generation of adaptive control systems. \\

\noindent \textit{What Challenges Still Remain?}

Despite its remarkable progress, extremum seeking is far from a completed discipline. Many important theoretical and practical challenges remain open.

One continuing objective is the development of algorithms capable of converging more rapidly while preserving robustness against measurement noise, disturbances, and model uncertainties. As applications become increasingly complex, there is also growing interest in reducing the computational burden associated with high-dimensional optimization and distributed decision-making.

Another important frontier concerns optimization under safety and operational constraints. Real engineering systems cannot freely explore every possible operating condition. Future extremum-seeking algorithms must therefore learn efficiently while respecting physical limitations, safety requirements, and communication constraints throughout the optimization process.

Finally, extending rigorous theoretical guarantees to increasingly sophisticated learning architectures remains a major challenge. Integrating extremum seeking with machine learning, event-triggered control, distributed optimization, stochastic systems, and large-scale cyber-physical infrastructures will require new mathematical tools that preserve the reliability that has always distinguished feedback control theory.

The next century of extremum seeking will likely be defined not by replacing classical control principles, but by expanding them into increasingly autonomous and interconnected engineering systems.

\section{Personal Reflections and Conclusions}

\noindent \textit{Why Has Extremum Seeking Endured for More Than a Century?}

Looking back over more than one hundred years of development, one cannot help but ask a simple question: \textbf{why has extremum seeking remained relevant for so long?} Few control methodologies have experienced such sustained growth, continually attracting new researchers while expanding into increasingly diverse application domains. 

The answer, we believe, lies in the remarkable permanence of the problem itself. 

Technology has changed dramatically since LeBlanc proposed the first extremum-seeking algorithm in 1922. Mechanical systems have given way to cyber-physical systems. Analog electronics have become digital computation. Today's controllers interact with cloud computing, autonomous vehicles, artificial intelligence, and massive communication networks. Yet despite these extraordinary technological transformations, engineers continue to face exactly the same fundamental question:
\begin{center}
\textit{How can a system improve its performance when the optimum is unknown?}
\end{center}
This question has never disappeared. On the contrary, it has become even more important as engineering systems have grown more complex.

Every generation has answered it differently. The pioneers relied primarily on engineering intuition. Later generations developed the mathematical tools necessary to establish rigorous convergence guarantees. Today, researchers are extending those same principles to distributed systems, infinite-dimensional dynamics, autonomous decision-making, and learning architectures.

The questions evolve. The mathematics evolves. The applications evolve. The central challenge does not.

Perhaps this explains why extremum seeking has repeatedly reinvented itself throughout its history. Rather than being tied to a particular mathematical technique or technological era, it addresses one of the most universal problems in engineering: learning how to do better while the system is operating. \\

\noindent \textit{Toward the Fourth Revolution.}

If the past century can be viewed as three successive scientific revolutions---the engineering revolution, the mathematical revolution, and the expansion toward increasingly complex dynamical systems---what might the next revolution look like?

We believe it will not simply consist of applying extremum seeking to new applications. Instead, it will involve a deeper integration of optimization, learning, prediction, communication, and decision-making into a unified feedback architecture.

Future autonomous systems will no longer perform isolated tasks such as regulation or optimization. They will continuously perceive their environment, predict future behavior, cooperate with other agents, adapt to unforeseen circumstances, and optimize multiple objectives simultaneously. In such systems, learning will not be an additional module attached to the controller; it will become an intrinsic component of the feedback loop itself.

Extremum seeking is uniquely positioned to contribute to this evolution. Unlike many modern learning paradigms, it was conceived from the outset as a feedback algorithm. Its language is that of stability, robustness, adaptation, and dynamical systems. These principles remain indispensable whenever autonomous decisions must be made safely in real time.

Rather than being replaced by artificial intelligence, extremum seeking may become one of the mechanisms through which intelligent systems acquire mathematical reliability. \\

\noindent \textit{A Personal Reflection.}

When I first became interested in extremum seeking, I was fascinated by the elegance of its central idea. A controller deliberately introduces a small perturbation, observes the system's response, and gradually discovers how to improve its own performance. Such a simple mechanism is capable of solving optimization problems that would otherwise appear remarkably difficult. 

Over the years, however, I came to realize that the true beauty of extremum seeking lies elsewhere.

It is not merely an optimization algorithm.

It is a philosophy of interaction with uncertainty.

Instead of insisting on complete knowledge before acting, extremum seeking accepts that uncertainty is inevitable. It experiments carefully, learns continuously, and improves incrementally. In doing so, it reflects a broader scientific principle: understanding is often achieved not by waiting for perfect models, but by intelligently interacting with the world.

Perhaps this is why the methodology has continued to evolve for more than a century. Its core idea is not tied to a specific application, mathematical technique, or technological generation. It captures a timeless way of solving problems---through curiosity, feedback, and adaptation. \\

\noindent \textit{Concluding Remarks.}

One hundred years ago, extremum seeking was introduced as an ingenious engineering solution to improve the operation of electrical systems. Few could have imagined that the same fundamental principle would eventually influence fields as diverse as robotics, renewable energy, process control, biomedical engineering, autonomous vehicles, distributed optimization, game theory, and infinite-dimensional systems.

Its history illustrates something that extends well beyond extremum seeking itself. Lasting scientific ideas rarely survive because they provide definitive answers. They survive because they continue to illuminate important questions.

The question that motivated LeBlanc's pioneering work in 1922 remains just as relevant today:
\begin{center}
\textit{How can a dynamical system learn to perform better without first knowing what ``better'' is?}
\end{center}
For more than a century, extremum seeking has provided one compelling answer.

We suspect that, in the century to come, it will inspire many more.

\appendix
\section{Appendix}

\subsection{Fundamentals of Extremum Seeking} \label{fundamentals}

This section provides an overview of the basic gradient-based version
of ESC with periodic signals for static maps (free of infinite-dimensional dynamics).

\begin{figure}
        \centering 
        \includegraphics[width=0.4\textwidth]{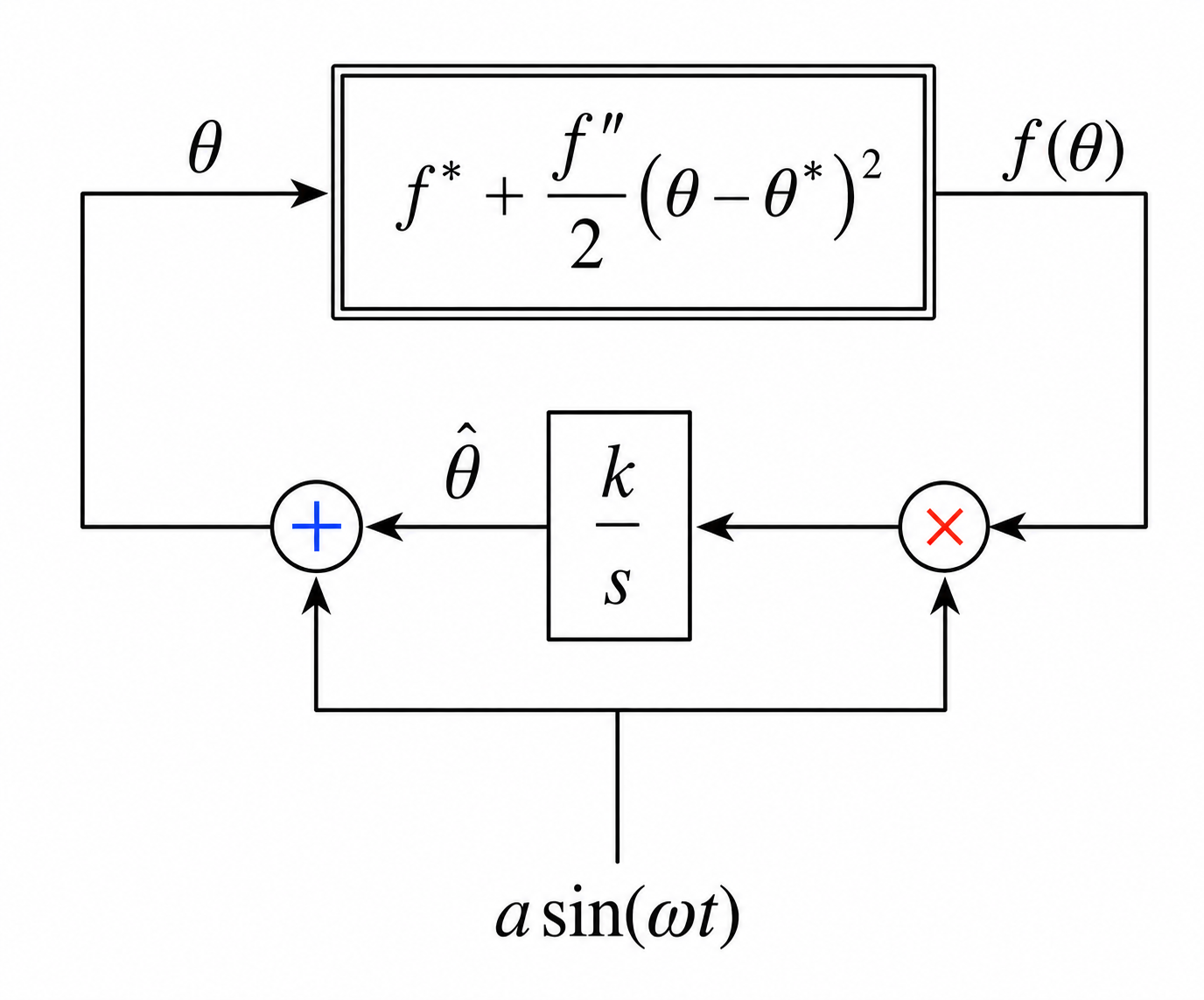}
        \caption{The simplest perturbation-based extremum seeking scheme for a quadratic single-input
map $f(\theta)$ in (\ref{Intro.quadratic_map}). The user has to only
know the sign of $f^{''}$, namely, whether the quadratic map has a maximum or a minimum, and
has to choose the adaptation gain $k$ such that $\sgn(k) = -\sgn\left(f^{''}\right)$. The user has to also choose
the frequency $\omega$ as relatively large compared to $a$, $k$ and $f^{''}$. Adapted from \cite{GKN:2012}.}
        \label{fig:controlador_basico_es}
    \end{figure}

Many versions of ESC exist, with various approaches to their study of stability 
\cite{KrsticWang2000}, \cite{LK:2012}, \cite{TNM:06}.
The most common version employs perturbation signals for the purpose of estimating
the gradient of the unknown map that is being optimized. To understand the basic idea of
extremum seeking, it is best to first consider the case of a static single-input map of the
quadratic form
\begin{equation} \label{Intro.quadratic_map}
f(\theta) = f^* + \frac{f^{''}}{2}(\theta - \theta^*)^2,
\end{equation}
where $f^*$, $f^{''}$ and $\theta^*$ are all unknown, as shown in Figure~\ref{fig:controlador_basico_es}. 

Three different $\theta$'s appear in Figure~\ref{fig:controlador_basico_es}:
$\theta^*$ is the unknown optimizer of the map, $\hat{\theta}(t)$
is the real-time estimate of $\theta^*$, and $\theta(t)$ is the actual input into the map. The actual input
$\theta(t)$ is based on the estimate $\hat{\theta}(t)$ but is perturbed by the signal $a \sin(\omega t)$ for the purpose
of estimating the unknown gradient $f^{''} \cdot (\theta - \theta^*)$ of the map $f(\theta)$ in (\ref{Intro.quadratic_map}).
The sinusoid is only one
choice for a perturbation signal---many other perturbations, from square waves to stochastic
noise, can be used in lieu of sinusoids, provided that they are of zero mean. The estimate $\hat{\theta} (t)$
is generated with the integrator $k/s$ with the adaptation gain $k$ controlling the speed of estimation.

The ESC algorithm is successful if the error between the estimate $\hat{\theta}(t)$ and the unknown $\theta^*$, namely the signal
    \begin{equation}
        \Tilde{\theta} (t) = \hat{\theta} (t) - \theta^*,
        \label{eq:erro_de_estimacao}
    \end{equation}
converges towards zero or some small neighborhood of zero as $t \!\to\! +\infty$. Based on Figure~\ref{fig:controlador_basico_es}, the estimate $\hat{\theta}(t)$ is governed by the differential equation $\dot{\hat{\theta}} = k \, \sin (\omega t) \, f(\theta)$, which means that the estimation error is governed by
    \begin{equation}
        \frac{d \tilde{\theta}}{dt} = 
            k \, a \, \sin (\omega t) 
            \bigg[ f^* + \frac{f^{''}}{2} 
                \big(\tilde{\theta} + a \, \sin (\omega t)
                \big)^2
            \bigg]\,.
        \label{eq:derivada_do_erro_de_estimacao}
    \end{equation}

Expanding the right-hand side, one obtains
    \begin{equation}\label{eq:slowandfastterms}
    \begin{split}
        \frac{d \tilde{\theta} (t)}{dt} =& 
            k \, a \, f^* \, \underbrace{\sin (\omega t)}_{\mbox{mean = 0}} +
            k \, a^3 \, \frac{f^{''}}{2} \underbrace{\sin^3 (\omega t)}_{\mbox{mean = 0}} \\
        &+  k \, a \, \frac{f^{''}}{2} 
                        \underbrace{ \sin (\omega t)}_{\mbox{fast, mean = 0}} \,\,\,\,\,\,\, \underbrace{\tilde{\theta} (t)^2}_{\mbox{slow}} \\
        &+  k \, a^2 \, f^{''} 
                        \underbrace{\sin^2 (\omega t)}_{\mbox{fast, mean = 1/2}} \,\,\,\,\,\,\, \underbrace{\tilde{\theta}(t)}_{\mbox{slow}}.
    \end{split}
    \end{equation}  
A theoretically rigorous time-averaging procedure \cite[Section 10.4]{khalil1996nonlinear} allows one to replace the above sinusoidal signals
by their means, yielding the ``average system''	\cite[p. 404]{khalil1996nonlinear}:
    \begin{equation}
        \frac{d \tilde{\theta}_{\rm{av}}}{dt} =
            \frac   {
                    \overbrace{k \, f^{''}}^{<0} \, 
                    a^2}
                    {2} \tilde{\theta}_{\rm{av}},
    \end{equation}
which is exponentially stable. The averaging theory guarantees that there exists a sufficiently
large $\omega$ such that, if the initial estimate $\hat{\theta}(0)$ is sufficiently close to the unknown $\theta^*$, one has
    \begin{equation}
        | \theta (t) - \theta^* | \leq 
        | \theta (0) - \theta^* | \, 
            e^{\frac{k f^{''} a^2}{2}t} + 
            \mathcal{O} \Bigg( \frac{1}{\omega} \Bigg) + \mathcal{O}(a)\,, \quad \forall t \geq 0\,.
        \label{eq:desigualdade_erro_estiamdo}
    \end{equation}
For the user, the inequality (\ref{eq:desigualdade_erro_estiamdo}) guarantees that, if $a$ is chosen small and $\omega$ is chosen
large, the input $\theta (t)$ exponentially converges to a small interval---of order $\mathcal{O}(\frac{1}{\omega} +a)$---around the unknown $\theta^*$ and,
consequently, the output $f(\theta(t))$  converges to the vicinity of the optimal output $f^*$.		

\newpage
\subsection{Predictor Feedback and Boundary Control of PDEs} \label{SBA_Directory}

\subsubsection{Basic Idea of Predictor Feedback Design with Actuator Delay} \label{predictor_appendix}

Here we consider the linear infinite-dimensional system
\begin{equation} \label{eq:2.1}
    \dot{X}(t)=AX(t)+BU(t-D),
\end{equation}
where $X \in \mathbb{R}^n$ $(A,\;B)$ is a controllable pair, and the input signal $U(t)$ is delayed by $D$ units of time.

 Given a stabilizing gain vector $K$ for the undelayed system, namely, given a
 vector $K$ such that the matrix $A+BK$ is Hurwitz, our wish is to have a control
 that achieves
 \begin{equation} \label{eq:2.2}
     U(t-D)=KX(t).
 \end{equation}
This control can be alternatively written as:
\begin{equation} \label{eq:2.3}
    U(t)=KX(t+D),
\end{equation}
and it appears to be nonimplementable since it requires future values of the state. However, with the variation-of-constants formula, treating the current state $X(t)$ as  the initial condition, we have
 \begin{equation} \label{eq:2.4}
     X(t+D)= e^{AD}X(t) + \int ^t_{t-D}e^{A(t-\theta)}BU(\theta)d\theta,\;\; \forall t \geq 0.
 \end{equation}
  This yields a feedback law
  \begin{equation} \label{eq:2.5}
      U(t) = K \Bigg[ e^{AD}X(t) + \int ^t_{t-D}e^{A(t- \theta)}BU(\theta)d\theta  \Bigg],\;\; \forall t \geq 0, 
  \end{equation}
 which is implementable, but it is infinite-dimensional, since it contains the distributed delay term involving past controls, $ \int_{t-D}^{t} e^{A(t-\theta)} BU(\theta)d\theta$. 
The closed-loop system is delay-compensated,
\begin{equation} \label{eq:2.6}
\dot{X}(t) = (A + BK)X(t), \quad t \geq D,
\end{equation}
but this is true only after the control “kicks in” at $t = D$. During the interval $t \in [0, D]$, the system state is governed by
\begin{equation} \label{eq:2.7}
X(t) = e^{At} X(0) + \int_{0}^{t} e^{A(t - \tau)} BU(\tau - D) \, d\tau, 
\quad \forall t \in [0,\; D].
\end{equation}
The feedback law (\ref{eq:2.5}) was introduced within the framework of “finite spectrum assignment” \cite{A121,A135} and the “reduction approach” \cite{a14}. In the Section~\ref{backstepping_appendix} we derive the same control law, but in a considerably more complicated way, which will pay dividends by providing us with an explicit Lyapunov–Krasovskii function and the ability to conduct stability analysis in the time domain, as discussed in \cite{K:2009}.

\subsubsection{Backstepping Boundary Control of Transport PDEs} \label{backstepping_appendix}

\begin{figure}
    \centering
    \includegraphics[width=0.7\linewidth]{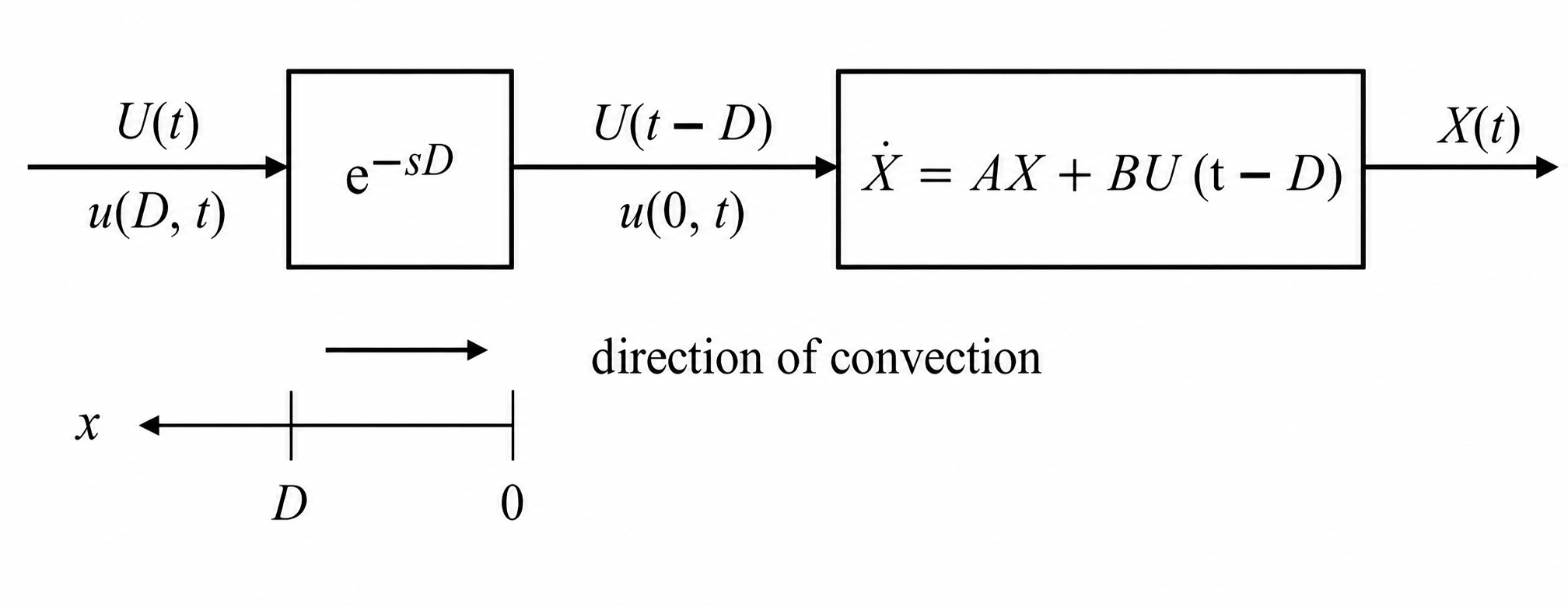}
    \caption{\textcolor{black}{Linear system $\dot{X}(t)=AX(t)+BU(t-D)$ with actuator delay $D$. Adapted from \cite[page 19]{K:2009}.}}
    \label{Fig2.1}
\end{figure}

The delay in the system (\ref{eq:2.1}) of the Section~\ref{predictor_appendix} can be modeled by the following first-order hyperbolic PDE, also referred to as the “transport PDE”:
\begin{eqnarray}  \label{eq:2.8}
    u_t(x,t) &=& u_x(x,t), \\
    u(D,t)&=&U(t). \label{eq:2.9}
\end{eqnarray}
The solution to this equation is
\begin{equation} \label{eq:2.10}
    u(x,t) = U(t + x - D), 
\end{equation}
and therefore the output
\begin{equation} \label{eq:2.11}
    u(0,t) = U(t - D) 
\end{equation}
gives the delayed input. The system (\ref{eq:2.1}) can now be written as
\begin{equation} \label{eq:2.12}
    \dot{X}(t) = AX(t) + Bu(0,t). 
\end{equation}
Equations \eqref{eq:2.8}–\eqref{eq:2.12} form an ODE–PDE cascade that is driven by the input $U$ from the boundary of the PDE (Figure~\ref{Fig2.1}).

\textcolor{black}{
Suppose a static state feedback control has been designed for a system with no delay (\textit{i.e.}, with $D = 0$) such that
\begin{equation}  \label{eq:eq:2.13}
    U(t) = KX(t) 
\end{equation}
is a stabilizing controller; \textit{i.e.}, the matrix $(A + BK)$ is Hurwitz. Consider the backstepping transformation
\begin{equation} \label{eq:2.14}
    w(x,t) = u(x,t) - \int_{0}^{x} q(x,y)u(y,t)\,dy - \gamma(x)^T X(t) 
\end{equation}
with which we want to map the system \eqref{eq:2.8}–\eqref{eq:2.12} into the target system
\begin{align} \label{eq:2.15}
    \dot{X}(t) &= (A + BK)X(t) + Bw(0,t), \\
    w_t(x,t) &= w_x(x,t), \label{eq:2.16} \\
    w(D,t) &= 0, \label{eq:2.17}
\end{align}
}
\textcolor{black}{
from the ODE-PDE system
\begin{align}
    \dot{X}(t) &= AX(t) + Bu(0,t), \label{eq:2.18} \\
    u_t(x,t) &= u_x(x,t), \label{eq:2.19} \\
    u(D,t) &= U(t). \label{eq:2.20}
\end{align}
The transformation \eqref{eq:2.14} applied to (\ref{eq:2.15})--(\ref{eq:2.20}) yields 
\begin{align}
    q_x(x,y) + q_y(x,y) &= 0, \label{eq:2.30} \\
    q(x,0) &= \gamma(x)^T B, \label{eq:2.31} \\
    \gamma'(x) &= A^T \gamma(x) \label{eq:2.32}
\end{align}
and 
\begin{equation}
    \gamma(0) = K^T. \label{eq:2.35}
\end{equation}
The solution to the ODE (\ref{eq:2.32}) is 
\begin{equation}
    \gamma(x)^T = K e^{Ax}. \label{eq:2.36}
\end{equation}
A general solution to \eqref{eq:2.30} is 
\begin{equation}
    q(x,y) = \phi(x - y), \label{eq:2.37}
\end{equation}
where the function $\phi$ is determined from \eqref{eq:2.31}. We get 
\begin{equation}
    q(x,y) = K e^{A(x - y)} B. \label{eq:2.38}
\end{equation}
We can now plug the gains $\gamma(x)$ and $q(x,y)$ into the transformation \eqref{eq:2.14} and set $x = D$ to get the control law:
\begin{equation}
    u(D,t) = \int_0^D K e^{A(D - y)} B u(y,t)\,dy + K e^{AD} X(t). \label{eq:2.39}
\end{equation}
The controller (\ref{eq:2.39}) is given in terms of the transport delay state $u(y,t)$. Using (\ref{eq:2.10}), one can also derive the representation in terms of the input signal $U(t)$:  
\begin{equation}
    U(t) = K \left[ e^{AD} X(t) + \int_{t-D}^t e^{A(t - \theta)} B U(\theta)\,d\theta\right], \label{eq:2.40}
\end{equation}
which is identical to the controller (\ref{eq:2.5}) in the Section~\ref{predictor_appendix}. }

\newpage
\subsection{Basic Idea of Nash Equlibrium Seeking in a Two-Player Game} \label{games_appendix}

Let players P1 and P2 represent two firms that produce the same good, have dominant control over a market, and compete for profit by setting their prices $u_1$ and $u_2$, respectively. The profit of each firm is the product of the number of units sold and the profit per unit, which is the difference between the sale price and the marginal or manufacturing cost of the product. In mathematical terms, the profits are modeled by
\begin{equation} \label{Intro_z.eq:1}
    J_i(t)=s_i(t)\big(u_i(t)-m_i \big),
\end{equation}
where $s_i$ is the number of sales, $m_i$ the marginal cost, and $i \in \{1,2\}$ for P1 and P2. Intuitively, the profit of each firm will be low if it either sets the price very low, since the profit per unit sold will be low, or if it sets the price too high, since then consumers will buy the other firm's product. The maximum profit is to be expected to lie somewhere in the middle of the price range, and it crucially depends on the price level set by the other firm.

To model the market behavior, we assume a simple, but quite realistic model, where for whatever reason, the consumer prefers the product of P1, but is willing to buy the product of P2 if its price $u_2$ is sufficiently lower than the price $u_1$. Hence, we model the sales for each firm as
\begin{equation} \label{Intro_z.eq:2}
    s_1(t)=S_d-s_2(t), \quad s_2(t)=\frac{1}{p}\big ( u_1(t)-u_2(t) \big ), 
\end{equation}
where the total consumer demand $S_d$ is held fixed for simplicity, the preference of the consumer for P1 is quantified by $p>0$, and the inequalities $u_1>u_2$ and $(u_1-u_2)/p < S_d$ are assumed to hold.

Substituting (\ref{Intro_z.eq:2}) into (\ref{Intro_z.eq:1}) yields expressions for the profits $J_1(u_1,u_2)$ and $J_2(u_1,u_2)$ that are both quadratic functions of the prices $u_1$ and $u_2$, namely,
\begin{eqnarray} \label{Intro_z.eq:4}
    J_1&=&\frac{-u_1^2+u_1u_2+(m_1+S_dp)u_1-m_1u_2-S_dpm_1}{p}, \\
    J_2&=&\frac{-u_2^2+u_1u_2+m_2u_1+m_2u_2}{p}, \label{Intro_z.eq:5}
\end{eqnarray}
and thus, the Nash equilibrium is easily determined to be
\begin{eqnarray} \label{Intro_z.eq:6}
    u_1^\ast&=&\frac{1}{3}(2m_1+m_2+2S_dp)\,, \\
    u_2^\ast&=&\frac{1}{3}(m_1+2m_2+S_dp)\,. \label{Intro_z.eq:7}
\end{eqnarray}
To make sure the constraints $u_1>u_2$, $(u_1-u_2)/p < S_d$ are satisfied by the Nash equilibrium, we assume that $m_1-m_2$ lies
in the interval $(-S_dp,2S_dp)$. If $m_1=m_2$, this condition is
automatically satisfied.

For completeness, we provide here the definition of a Nash
equilibrium $u^\ast=[u_1^\ast\,, u_2^\ast]^T$ in a 2-player game:
\begin{equation} \label{Intro_z.eq:8}
    J_i(u_i^\ast,u_{-i}^\ast) \ge J_i(u_i,u_{-i}^\ast)\,, \qquad \forall u_i \in U_i,\; i \in \{1, 2\}\,,
\end{equation}
where $J_i$ is the payoff function of player $i$, $u_i$ its action, $U_i$ its action set, and $u_{-i}$ denotes the action of the other player. Hence, no player has an incentive to unilaterally deviate its action from $u^\ast$. 
In the duopoly example, $U_1=U_2=\mathbb{R}_+$, where $\mathbb{R}_+$ denotes the set of positive real numbers.
\medskip

    \begin{figure}
        \centering
        \includegraphics[width=0.3\textwidth]{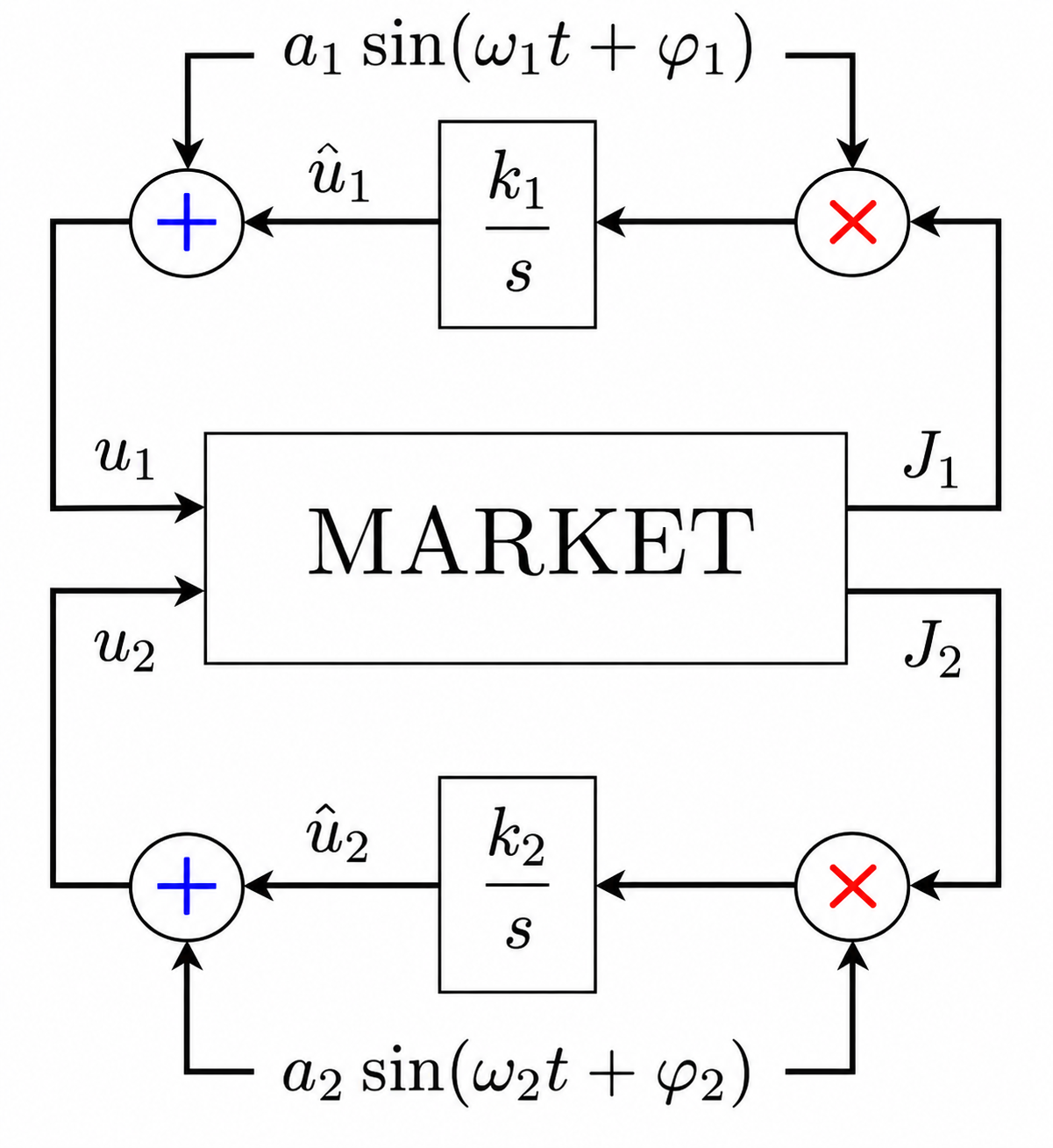}
        \caption{Deterministic Nash seeking schemes applied by players in a duopoly market structure. Adapted from \cite{FKB:2012}.}
        \label{fig_intro_z_Nash_Market}
    \end{figure}

To attain the Nash strategies (\ref{Intro_z.eq:6})--(\ref{Intro_z.eq:7}) without any knowledge of modeling information, such as the consumer's preference $p$, the total demand $S_d$, or the other firm's marginal cost or price, the firms implement a non-model based real-time optimization strategy, \textit{e.g.}, deterministic extremum seeking with sinusoidal perturbations, to set their price levels. Specifically, P1 and P2 set their prices, $u_1$ and $u_2$ respectively, according to the time-varying strategy \textcolor{black}{(Figure~\ref{fig_intro_z_Nash_Market}):}
\begin{eqnarray} \label{Intro_c.eq:00}
    \dot{\hat{u}}_i(t)&=&k_i\mu_i(t)J_i(t), \\
    u_i(t)&=&\hat{u}_i(t)+\mu_i(t),\label{Intro_c.eq:01}
\end{eqnarray}
where $\mu_i(t)=a_i\sin(\omega _i t + \varphi_i)$, $k_i$, $a_i$, $\omega_i > 0$ and $i \in \{1,2 \}$. 

\newpage
Further, the frequencies are of the form
\begin{equation} \label{Intro_c.eq:02}
    \omega_i = \omega \overline{\omega}_i,
\end{equation}
where $\omega$ is a positive real number and $\overline{\omega}_i$ is a positive rational number. 

In contrast, the firms are also guaranteed to converge to the Nash equilibrium when employing the standard parallel action update scheme \cite[Proposition $4.1$]{Basar:1999}
\begin{eqnarray} \label{Intro_z.eq:12}
    u_1^{(k+1)}&=&\frac{1}{2}\Big(u_2^{(k)}+m_1+S_dp \Big), \\
     u_2^{(k+1)}&=&\frac{1}{2}\Big(u_1^{(k)}+m_2 \Big), \label{Intro_z.eq:13}
\end{eqnarray}
which requires each firm to know both its own marginal cost and the other firm's price at the previous step of the iteration, and also requires P1 to know the total demand $S_d$ and the consumer preference parameter $p$. In essence, P1 must know nearly all the relevant modeling information. When using the extremum seeking algorithm (\ref{Intro_c.eq:00})--(\ref{Intro_c.eq:01}), the firms only need to measure the value of their own payoff functions, $J_1$ and $J_2$. Convergence
of (\ref{Intro_z.eq:12})--(\ref{Intro_z.eq:13}) is global, whereas the convergence of the Nash seeking strategy for this example can be proved to be semi-global, following \cite{TNM:06}, or locally, by applying the theory of averaging \cite{khalil1996nonlinear}.

\printcredits

\newpage
\section*{Acknowledgments}




The author is grateful to Prof. Miroslav Krstic, whose rare ability to recognize the future hidden inside a simple idea has been a constant \textit{source} of inspiration. More often than not, that vision has made the difference between an interesting discussion and an entirely new research direction. 
%

This work was supported in part by the Brazilian funding agencies Conselho Nacional de Desenvolvimento Científico e Tecnológico (CNPq) and Fundação Carlos Chagas Filho de Amparo à Pesquisa do Estado do Rio de Janeiro (FAPERJ). This study was also financed in part by the Coordenação de Aperfeiçoamento de Pessoal de Nível Superior – Brasil (CAPES) – Finance Code 001.

\section*{CRediT Author Statement}

\textbf{Tiago Roux Oliveira} is the sole author of this manuscript and was responsible for its
conceptualization, literature investigation, development of the historical and technical
perspective, preparation and adaptation of the illustrative material, writing of the original
draft, and review and editing of the final manuscript.

\bibliographystyle{elsarticle-num}

\bibliography{cas-refs}

\newpage
\vskip6pc         

\bio{figs/TRO}
\textbf{Tiago Roux Oliveira} (Senior Member, IEEE) received the B.Sc. degree in Electrical Engineering from the Rio de Janeiro State University (UERJ), Rio de Janeiro, Brazil, in 2004, and the M.Sc. and D.Sc. degrees in Electrical Engineering from the Federal University of Rio de Janeiro (COPPE/UFRJ), Rio de Janeiro, Brazil, in 2006 and 2010, respectively. In 2014, he was a Visiting Scholar with the University of California, San Diego (UC San Diego), La Jolla, CA, USA. He is currently a Professor with the Department of Electronics and Telecommunications Engineering, UERJ. He has authored or coauthored approximately 300 peer-reviewed journal articles, conference papers, and book chapters. 
Dr. Oliveira is an Associate Editor of the Journal of the Franklin Institute, Systems $\&$ Control Letters, International Journal of Robust and Nonlinear Control, IEEE Open Journal of Control Systems, IEEE Transactions on Cybernetics, IEEE Control Systems Letters, IEEE Transactions on Automatic Control, and Automatica. 
His research interests include adaptive and learning control, extremum seeking, distributed optimization, game-theoretic learning, and control of systems with delays and partial differential equations.
His research contributions have been recognized with several awards, including the Brazilian CAPES National Best Ph.D. Thesis Award in Electrical Engineering (2011), the Young Scientist of Our State Award (JCNE) from FAPERJ (2012, 2015, and 2018), the Scientist of Our State Award (CNE) from FAPERJ (2021 and 2024), and the 2021 IEEE Transactions on Control Systems Technology Outstanding Paper Award. 
He has served on several IFAC Technical Committees, including Adaptive and Learning Systems (TC 1.2), Control Design (TC 2.1), Nonlinear Control Systems (TC 2.3), and Distributed Parameter Systems (TC 2.6), as well as on the IEEE Control Systems Society Technical Committees on Variable Structure and Sliding Mode Control and Distributed Parameter Systems. In 2017, he was elected an Affiliated Member of the Brazilian Academy of Sciences. He served as Chair of the IFAC Technical Committee on Adaptive and Learning Systems (TC 1.2) from 2020 to 2026 and as President of the Brazilian Society of Automatics (SBA) from 2023 to 2025. 
He is the coauthor, with Miroslav Krstic, of the book Extremum Seeking through Delays and PDEs (SIAM, 2022), which presents a unified treatment of extremum seeking for systems with delays, partial differential equations, and distributed optimization.
\endbio

\end{document}